\documentclass[ejs]{imsart}

\usepackage{amsthm,amsmath}
\RequirePackage[numbers]{natbib}
\RequirePackage[colorlinks,citecolor=blue,urlcolor=blue]{hyperref}

\usepackage{bbm}
\usepackage{amssymb}
\usepackage{graphicx}
\usepackage{url}
\doi{00000000}
\pubyear{0000}
\volume{0}
\firstpage{0}
\lastpage{0}

\startlocaldefs

\def\sig{\sigma}
\def\Sig{\Sigma}
\def\gam{\gamma}

\newtheorem{condition}{Condition}[section]
\newtheorem{theorem}{Theorem}[section]

\newtheorem{lemma}[theorem]{Lemma}

\newcommand{\bel}{\begin{eqnarray}\label}
\newcommand{\eel}{\end{eqnarray}}
\newcommand{\bes}{\begin{eqnarray*}}
\newcommand{\ees}{\end{eqnarray*}}
\newcommand{\bei}{\begin{itemize}}
\newcommand{\beiftnt}{\begin{itemize}\footnotesize}
\newcommand{\eei}{\end{itemize}}

\def\benu{\begin{enumerate}}
\def\eenu{\end{enumerate}}

\def\argmin{\mathop{\rm arg\, min}}
\def\real{{\mathbb{R}}}
\def\R{{\real}}

\def\E{{\mathbb{E}}}
\def\P{{\mathbb{P}}}

\def\complex{\mathop{{\rm I}\kern-.58em\hbox{\rm C}}\nolimits}

\def\sgn{\hbox{\rm sgn}}

\def\Var{\hbox{\rm Var}}

\def\eps{\epsilon}

\def\lam{\lambda}

\endlocaldefs

\begin{document}

\begin{frontmatter}

\title{Debiasing the Debiased Lasso with bootstrap }
\runtitle{Bootstrapping the Debiased Lasso}


\author{\fnms{Sai} \snm{Li}\corref{}\ead[label=e1]{saili@pennmedicine.upenn.edu}}
\address{Department of Biostatistics, University of Pennsylvania, Philadelphia, PA19104 }

\runauthor{S.Li}

\begin{abstract}
We consider statistical inference for a single coordinate of regression coefficients in high-dimensional linear models. 
Recently, the debiased estimators are popularly used for constructing confidence intervals and hypothesis testing in high-dimensional models. However, some representative numerical experiments show that they tend to be biased for large coefficients, especially when the number of large coefficients dominates the number of small coefficients. In this paper, we propose a modified debiased Lasso estimator based on bootstrap. Let us denote the proposed estimator BS-DB for short. We show that, under the irrepresentable condition and other mild technical conditions, the BS-DB has smaller order of bias than the debiased Lasso in existence of a large proportion of strong signals. If the irrepresentable condition does not hold, the BS-DB is guaranteed to perform no worse than the debiased Lasso asymptotically. Confidence intervals based on the BS-DB are proposed and proved to be asymptotically valid under mild conditions.
Our study on the inference problems integrates the properties of the Lasso on variable selection and estimation novelly.
The superior performance of the BS-DB over the debiased Lasso is demonstrated via extensive numerical studies.
   
\end{abstract}





\end{frontmatter}





\section{Introduction}

\subsection{Background}
High-dimensional linear models have broad applications in many fields, such as biology, genetics, and machine learning.  
A number of statistical methods have been introduced to solve the problems on prediction, estimation, and variable selection regarding regression coefficients. On the other hand, statistical inference in high-dimensional models has recently caught a lot of research interests and efforts for its importance in providing uncertainty assessment and the nontrivial statistical challenges.



The Lasso estimator \citep{Tibshirani96} has been a popular tool for modeling high-dimensional data. When the number of covariates $p$ is fixed, however, it has been shown to have no closed form for its limiting distribution in the low dimensional setting \citep{KF00}. \citet{CL10} showed the inconsistency of bootstrapping the Lasso if at least one coefficient is zero. Thus, there is substantial difficulty in drawing valid inference based on the Lasso estimates directly. 
 Nevertheless, \citet{CL11} developed a modified bootstrap estimator based on the Lasso as well as a bootstrap estimator based on  Adaptive Lasso \citep{Zou06}.
For $p$ increasing with the sample size $n$, \citet{CL13} showed the bootstrap approximation consistency for Adaptive Lasso estimators under some technical conditions. 

In the $p\gg n$ scenario, \citet{ZZ14} proposed asymptotically Gaussian-distributed estimators of low-dimensional parameters in high-dimensional linear regression models. Such estimators are known as the ``debiased Lasso'' or the ``de-sparsifying Lasso''. In the same paper, they proposed an optimization scheme for calculating the ``correction score'', whose properties are carefully studied in \citep{JM14} for both fixed designs and sub-Gaussian designs. Along this line of research, many recent papers have studied relevant generalizations for the debiased approach. 
\citet{Van14} considered the debiased Lasso estimator in generalized linear models with convex loss functions. 
\citet{Buh15} and \citet{JV15} studied statistical inference in misspecified high-dimensional linear models and graphical models, respectively.
\citet{FNL16} considered the debiased method in high-dimensional Cox models. From the minimax perspective, \citet{CG17} studied the optimal expected lengths of confidence intervals for linear combinations of regression coefficients in sparse high-dimensional linear models.
\citet{JM15} considered sample size conditions for the debiased Lasso method in high-dimensional linear models with Gaussian design and Gaussian noise. Under some regularity conditions, they show a potentially weaker sample size condition when the true precision matrix of the design is sparse. 
 Related approaches are also actively studied in the context of econometrics and causal inference \citep{BCH14, BCW16, Chern18}. 


The present work is motivated by the connections between statistical inference and variable selection problems. 
Variable selection has become an active research topic in high-dimensional literature for decades. Many established variable selection methods have been proposed and studied \citep{Tibshirani96, Candes07, FL01,Wasserman09, Meinshausen10, MCP}.
It is known that if the nonzero coefficients can be consistently selected, least square estimators based on the selected model can lead to asymptotically valid inference procedures. However, the consistency of variable selection always requires the beta-min condition, which assumes  the strengths of nonzero coefficients are uniformly larger than certain threshold. This condition is uncheckable and can be hard to fulfill in applications. In a recent paper \citep{Zhao17}, the post-Lasso least squares is justified for asymptotic valid statistical inference in high-dimensional linear models. Their analysis is based on the conditions guaranteeing the set of variables selected by the Lasso is deterministic with high probability. In the current work, we explore the interaction between variable selection and inference problems and demonstrate the benefits of having a large proportion of strong signals for statistical inference under proper conditions. Specially, we do not require the beta-min condition and the selected set of variables is not necessarily deterministic.

Another philosophy for inference in the high-dimensional setting is based on selective inference, whose focus is on making inference conditional on the selected model \citep{LTTT14,CHS15,Lee16,TTLT16}. However, it is not considered in the current work. 

Our proposed approach is closely related to the bootstrap procedures for inference. Bootstrap has been widely used in high-dimensional models for conducting statistical inference. \citet{Mam93} considered estimating the distribution of linear contrasts and of F-test statistics when $p$ increases with $n$. \citet{CCK13} developed theories for multiplier bootstrap to approximate the maximum of a sum of high-dimensional random vectors. \citet{Deze17} proposed residual, paired and wild multiplier bootstrap methodology for the debiased Lasso estimators. \citet{ZC17} proposed a bootstrap-assisted debiased Lasso estimator to conduct simultaneous inference for non-Gaussian errors. The purpose of using bootstrap in aforementioned two papers mainly concern dealing with heteroscedastic errors as well as simultaneous inference. In the present work, we show the bias correction effect of bootstrap in high-dimensional inference. 

\subsection{The debiased approach}
\label{sec1-3}
The debiased Lasso \citep{ZZ14, Van14, JM14}  for high-dimensional linear models can be described as follows.
Consider a linear regression model 
\begin{equation}
\label{lm1}
   y_i = x_i^T\beta+\epsilon_i,
\end{equation}
where $\beta \in \R^p$ is a vector of regression coefficients and $\epsilon_1,\dots,\epsilon_n$ are \textit{i.i.d.} random variables with mean 0 and variance $\sigma^2$. We consider the high-dimensional scenario where $p$ can be larger or much larger than $n$. We assume $\beta$ is sparse with support $S$ such that $|S|=s$. Let $X\in\R^{n\times p}$ be the design matrix with the $i$-th row being $x_i^T$ and $y=(y_1,\dots,y_n)^T$. Let $X_j$ denote the $j$-th column of $X$. Let $\Sig^n=X^TX/n$.  Let $\Sig=\E[\Sig^n]$ denote the population gram matrix which is positive definite and $\Theta=\Sig^{-1}$. Let $\Theta_j=\Theta_{.,j}$ denote the $j$-th column of $\Theta$. 

 The Lasso \citep{Tibshirani96} estimator of $\beta$ is defined as
 \begin{equation}
\label{beta-init}
  \hat{\beta}=\argmin_{b\in \R^p}\left\{ \frac{1}{2n}\|y-Xb\|_2^2+\lam\|b\|_1\right\}
\end{equation}
for some tuning parameter $\lam>0$. Some of its variations have been proposed and studied \citep{Zou06, squared-lasso, scaled-lasso}.
Suppose that we are interested in making inference of $\beta_j$ for some $j\in\{1,\dots, p\}$. 
 The debiased Lasso estimator of $\beta_j$ can be written as
\begin{equation}
\label{beta-db}
   \hat{\beta}_j^{(DB)}=\hat{\beta}_j+\frac{\langle \hat{z}_j,y-X\hat{\beta}\rangle}{n},
\end{equation}
where $\hat{z}_j\in\R^n$ is the so-called correction score and can be computed via another Lasso regression. 
Specifically, define 
\begin{align}
\hat{\gam} &= \argmin_{\gam\in \R^{p}} \left\{\frac{1}{2n}\|X\gamma\|_2^2+ \lam_j\|\gamma_{-j}\|_1:~\gamma_j=1\right\}\label{hgam0}
\end{align}
for some tuning parameter $\lam_j>0$ and 
\begin{equation}
\label{eq-z1}
  \hat{z}_{j}=\left\{\begin{array}{ll}
  X\Theta_j~\text{if $\Theta$ is known}\\
  \frac{nX\hat{\gam}}{\langle X_j,X\hat{\gam}\rangle}~\text{if $\Theta$ is unknown}.
   \end{array}
   \right.
  \end{equation}
In fact, the correction score can also be realized via a quadratic optimization approach \citep{ZZ14, JM14} when $\Theta$ is unknown. The difference is that the optimization in (\ref{hgam0}) can utilize the sparsity of $\Theta_j$, if $\Theta_j$ is sparse indeed, and can achieve semiparametric efficiency under proper conditions \citep{Van14}. 

A $100\times (1-\alpha)\%$ two-sided confidence interval for $\beta_j$ can be constructed as
\begin{align}
\label{ci0}
\hat{\beta}_j^{(DB)}\pm q_{\alpha/2}\hat{\sig}\|\hat{z}_j\|_2/n,
\end{align}
where $q_{\tau}$ is $\tau$-th quantile of standard normal distribution and $\hat{\sig}$ is some consistent estimate of the noise level $\sig$.

When $\Theta$ is unknown, it has been proved that the asymptotic normality of $\hat{\beta}_j^{(DB)}$ requires $n\gg (s\log p)^2$ and some other technical conditions, say, in \citet{Van14}. These conditions guarantee that the remaining bias of the debiased Lasso is asymptotically sufficiently small such that $n^{-1/2}$-length confidence intervals are achievable. In \citet{CG17}, it has been shown that the minimax optimal confidence intervals for $\beta_j$ has lengths of order $n^{-1/2}+s\log p/n$. That is, to achieve a confidence interval with length of order $n^{-1/2}$, the condition $n\gg (s\log p)^2$ is unavoidable in the minimax sense. This reveals the optimality of the debiased Lasso procedure. 

Actually, when investigating the worst case scenario considered in \citet{CG17} (Theorem 3), one can see that it concerns the case where $\underline{s}=\large|\{j: |\beta_j|\asymp \sqrt{\log p/n}\}\large|\asymp s$. That is, loosely speaking, the number of ``weak'' coefficients dominates the number of ``strong'' coefficients. Therefore, if there exists a large proportion of strong signals, the debiased Lasso may not be optimal. If not, it is unknown what is a better procedure for statistical inference. This question is of significant practical value in view of the following numerical experiments, which suggest that the debiased Lasso estimators can be severely biased for large signals.

\subsubsection{The bias of the debiased Lasso related to signal strengths}
We consider a simplified setting where the observations are generated from model (\ref{lm1}) with known noise level $\sig=1$ and known $\Sig$. We set $n=100$ and $p=300$. Each row of $X$ is \textit{i.i.d.} from a multivariate Gaussian distribution with mean zero and covariance matrix $\Sig$. The noise vector $\eps$ are \textit{i.i.d.} from standard Gaussian distribution. We consider five levels of sparsity, i.e. $s=2k$ for $k=2,\dots,8$. For each integer $2\leq k\leq 8$, we set $S=\{1,\dots,2k\}$ and $\beta_S=\left(4,2,4,2,\dots,4,2,4,0.2\right)^T$.  Comparing with the inflated noise level $\sig\sqrt{2\log p/n}\approx 0.338$, first $s-1$ coefficients can be viewed as strong signals and the $s$-th coefficient can be viewed as weak. We consider two formats of $\Sig$. The first one is $\Sig=I_{p}$, i.e. the identity matrix. The second one is $\Sig=\Sig^o$, where $\Sig^o$ is a block diagonal matrix with $\Sig^o_{S^c,S^c}=I_{p-s}$ and $\Sig^o_{S,S}$ is Toeplitz with the first row equals $(1,-0.1,\mathbf{0}_{s-3}^T, -0.1)$. 
We compute the debiased Lasso (DB) based on (\ref{beta-init}) and (\ref{beta-db}) with tuning parameter $\lam=2\sig\sqrt{\log p/n}$ and the oracle correction sore $\hat{z}_j=X\Theta_j$.

\begin{figure}
\centering
\includegraphics[height=5cm]{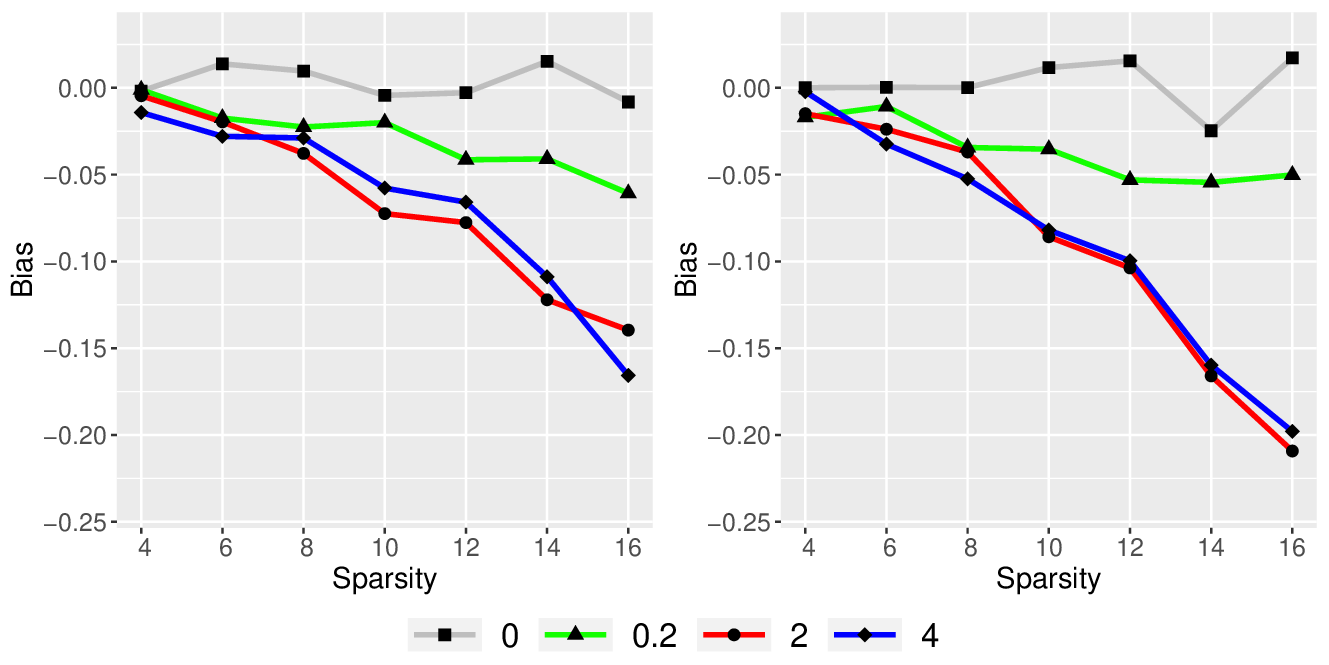}
\caption{Estimation errors of the debiased Lasso for $\beta_j\in\{4,2,0.2,0\}$ at different sparsity levels with $\Sig=I_p$ (left panel) and $\Sig=\Sig^o$ (right panel). Each point is the median estimate based on 500 independent experiments.}
\label{fig0}
\end{figure}
When $\Sig$ is identity (left panel of Figure \ref{fig0}), distinct patterns are observed for signals with different strengths. The debiased Lasso for $\beta_j$, $j\in S^c$ has bias floating around zero. For $\beta_j$, $j\in S$, the debiased estimators are always negatively biased across different sparsity levels. This implies that the debiased Lasso can lead to low chance of discovering true signals under the current set-up.
Moreover, the debiased estimators for strong signals are more severely biased as the sparsity level $s$ gets larger.  The right panel of Figure \ref{fig0} displays the bias of the debiased Lasso when $\Sig=\Sig^o$. The difference of the left and right panel of Figure \ref{fig0} preludes the effect of $\rho_S$ defined in (\ref{rho-s}) where $\rho_S=1$ with $\Sig=I_p$ and $\rho_S=1.25$ with $\Sig=\Sig^o$ at all sparsity levels. When $\rho_S$ gets larger, the bias of the debiased Lasso on the true support gets larger.

Motivated by this numerical experiment, we look into the effects of signal strengths on statistical inference in high-dimensional models. Under proper conditions, we provide a new error analysis of the debiased Lasso which shows its distinct behaviors on and off the true support. 

\begin{figure}
\centering
\makebox{
 \includegraphics[height=5cm]{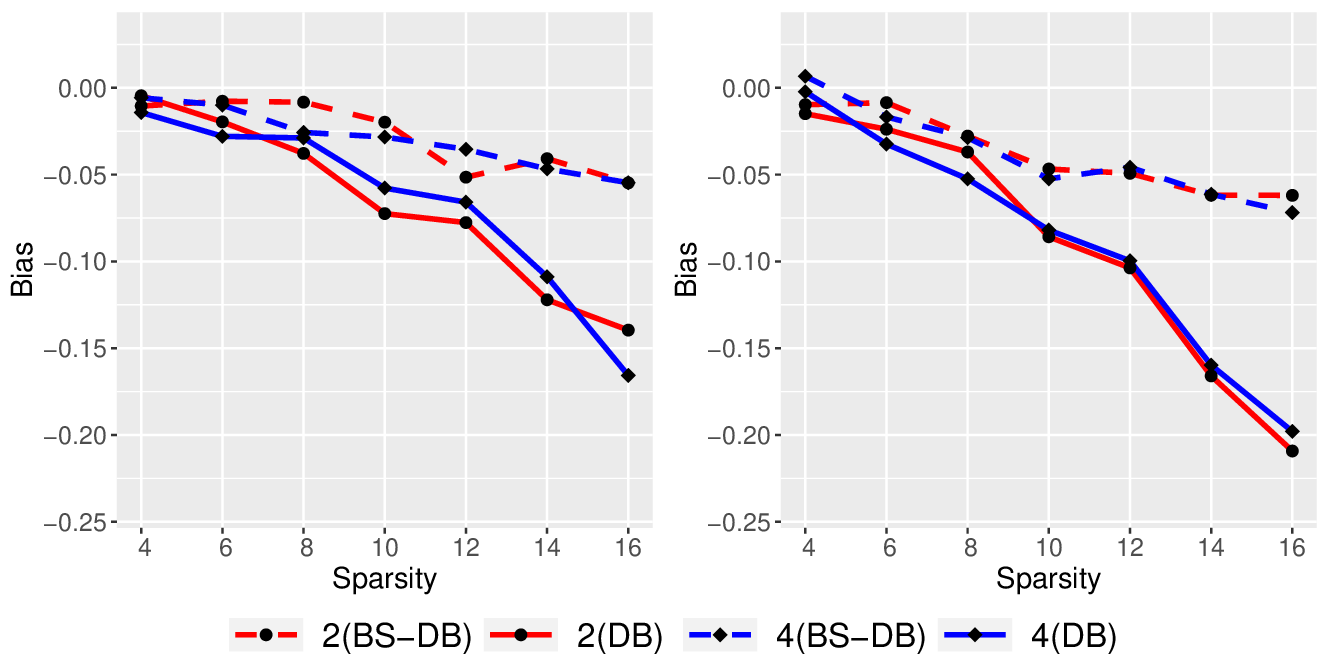}
}
\caption{Estimation errors of the debiased Lasso (DB) and of the BS-DB for $\beta_j\in\{4,2\}$ at different sparsity levels with identify covariance matrix (left panel) and $\Sig=\Sig^o$ (right panel). Each point is the median estimate based on 500 independent experiments.}
\label{fig01}
\end{figure}

More importantly, we introduce a bootstrapped debiased Lasso approach (BS-DB), which can have smaller order of bias than the debiased Lasso when there are a large proportion of strong signals. 
 Figure \ref{fig01} unveils the bias correction effect of the BS-DB estimator in comparison to the debiased Lasso in the simulation settings considered above. It is not hard to see that the bias for strong signals are significantly reduced with our proposal. The BS-DB estimator is constructed according to Section \ref{sec2-1} with the same tuning parameter as in debiased Lasso.
More numerical results on constructing confidence intervals with unknown covariance matrix and unknown noise level are presented in Section \ref{sec-simu}. In next subsection, we summarize the major contributions of the current work.

 \subsection{Summary of our contributions}
 For the proposed BS-DB estimator defined in (\ref{beta-bsdb}), we will prove that 
\begin{align}
\label{eq-bsdb0}
&\hat{\beta}^{(BS-DB)}_j-\beta_j=\widehat{rem}_j+\hat{\eta}_j,
\end{align}
where $\widehat{rem}_j$ is the remaining bias of $\hat{\beta}^{(BS-DB)}_j$ and $\sqrt{n}\hat{\eta}_j$ is asymptotically normal. We upper bound the magnitude of $\widehat{rem}_j$ assuming sub-Gaussian noise (Condition \ref{cond2b}). 

\textbf{(a)} Suppose that the design matrix is row-wise independent Gaussian (Condition \ref{cond1b}) and the irrepresentable condition (Condition \ref{cond3b}) holds, if $n\gg s\log p$, then
\begin{align}
\label{eq-re2}
|\widehat{rem}_j|=_{O_p}\left\{
\begin{array}{ll}
\frac{\sqrt{s_-\log p}}{n}+ \frac{\rho_S\sqrt{s_-\log p}}{n}\mathbbm{1}(j\in S)~~\text{if $\Theta$ is known}\\
\frac{s_-\log p}{n}~~\text{if $\Theta$ is unknown},
\end{array}\right.
\end{align}
where $\mathbbm{1}(\cdot)$ is an indicator function, $s_-=\left|\{j\in S:|\beta_j|\leq C\rho_S\lam\}\right|$ for a large enough constant $C$, and $\rho_S$ defined in (\ref{rho-s}). It implies that, under proper conditions, asymptotically valid confidence intervals can be constructed based on $\hat{\beta}^{(BS-DB)}_j$ if $s\log p\lesssim n$ and $(s_-\log p)^2\ll n$ for unknown $\Theta$. This sample size condition is weaker than the one for the debiased Lasso if $s_-\ll s$. 
We will bring more discussions on the demand and potential relaxations of the irrepresentable condition in Section \ref{sec3}.
As a byproduct, we get a new error expansion of the debiased Lasso under current conditions in Theorem \ref{thm2a}, which shows that its remaining bias can be of different magnitude for $j\in S$ and $j\notin S$.

\textbf{(b)} When the design is sub-Gaussian (Condition \ref{cond4b}), without assuming the irrepresentable condition, we prove that
 the remaining bias of the BS-DB is $O_P(s\log p/n)$ under mild conditions, which is is asymptotically ``no worse'' than the debiased Lasso (Theorem \ref{thm10}).

\subsection{Notations}
For a set $A\subseteq \{1,\dots,p\}$, $X_A\in \R^{n\times |A|}$ is the submatrix formed by $X_j,~j\in A$. Let $A^c=\{1,\dots,p\}\setminus A$. For $j=1,\dots,p$, let $\{-j\}=\{1,\dots,p\}\setminus \{j\}$.
For a vector $v\in \R^k$, let $\|v\|_p$ denote the standard $\ell_p$-norm of $v$. For another vector $u\in \R^k$, let $\langle v,u\rangle = \sum_{i=1}^k v_iu_i$. 
For a symmetric matrix $D\in \R^{k_1\times k_1}$, let $\Lambda_{\max}(D)$ and $\Lambda_{\min}(D)$ be the largest and smallest eigenvalues of $D$. Let $\textrm{Tr}(D)$ denote the trace of $D$.
Let $I_p$ be the $p\times p$ identity matrix.  We use $e_j$ to refer to the $j$-th standard basis element, i.e. $e_1 = (1,0,\dots,0)^T$. We use $e_{S_1}\in\R^{p\times |S_1|}$ to refer to the sub-matrix of an identify matrix, which contains $j$-th column of $I_{p}$ for $\forall j\in S_1$.
For a random variable $v$, let $\|v\|_{\psi_2}=\sup_{t\geq 1}t^{-1/2}\left(\E[|v|^t]\right)^{1/t}$ denote its sub-Gaussian norm. For a random vector $V\in\R^k$, let $V_{\psi_2}=\sup_{\|b\|_2=1}\|\langle V,b\rangle\|_{\psi_2}$.
Let $f(n)$ and $g(n)$ be two functions. We use $f(n)\ll g(n)$ to refer to ``$f(n)=o(g(n))$''. The notation $f(n)\gg g(n)$ can be defined analogously. We use $f(n)\lesssim g(n)$ to refer to $f(n)=O(g(n))$ and $f(n)\gtrsim g(n)$ can be defined analogously.
We use $c,c_0,c_1,\dots$ and $C_0,C_1,\dots$ to denote generic constants that can vary from one position to the other. 
\subsection{Organization of the rest of the paper}
The rest of the paper is organized as follows. In Section \ref{sec2}, we introduce the proposed approach for constructing confidence intervals in high-dimensional linear models. In Section \ref{sec3}, we prove the theoretical properties of the proposed approach with and without the irrepresentable condition. In Section \ref{sec-simu}, we demonstrate the empirical performance of the BS-DB for statistical inference in comparison to the debiased Lasso in various settings. In Section \ref{sec-diss}, we bring more discussions to the topics related to this paper. The proofs are provided in Section \ref{sec-pf} and in the Appendix. 

\section{Bootstrapping the debiased Lasso}
\label{sec2}
In this section, we introduce the proposed procedure and bring some intuitions to its merits. 

\subsection{The procedure for constructing confidence intervals}
\label{sec2-1}

\begin{enumerate}
\item[(i)] (Fitting the modified debiased Lasso)
Let $\hat{\beta}$ be the initial Lasso estimator computed via (\ref{beta-init}) with tuning parameter $\lam$.
Let $\hat{w}_j$ be the correction score defined in (\ref{eq-z2}) (with known or unknown precision matrix). We compute a modified debiased Lasso estimator based on $\hat{\beta}$ and $\hat{w}_j$:
\begin{align}
\label{beta-db2}
\hat{\beta}^{(mDB)}_j=\hat{\beta}_j+\frac{\langle \hat{w}_j,y-X\hat{\beta}\rangle}{n}.
\end{align}
\item[(ii)] (Bootstrapping the Lasso)
Let $\hat{\beta}^*$ be the Lasso estimate with input $(X,X\hat{\beta})$ and tuning parameter $\lam$.
We estimate the bias of $\hat{\beta}^{(mDB)}_j$ with
\begin{align}
\label{eq-hbj}
\hat{b}^{*}_j=\hat{\beta}^*_j+\frac{\langle \hat{w}_j,X(\hat{\beta}-\hat{\beta}^*)\rangle}{n}-\hat{\beta}_j.
\end{align}
Substracting $\hat{b}^*_j$ from $\hat{\beta}_j^{(mDB)}$, we arrive at a bootstrapped debiased Lasso estimator
\begin{align}
\label{beta-bsdb}
\hat{\beta}_j^{(BS-DB)}=\hat{\beta}^{(mDB)}_j-\hat{b}^*_j.
\end{align}

\item[(iii)] 
An $100\times (1-\alpha)\%$ two-sided confidence interval for $\beta_j$ can be constructed as
\begin{align}
\label{ci1}
\hat{\beta}^{(BS-DB)}_j\pm q_{\alpha/2}\hat{\sig}\|\hat{w}_{j}\|_2/n,
\end{align}
where $\hat{\sig}$ can be any consistent estimator of $\sig$.
\end{enumerate}


The modified debiased Lasso computed in Step (i) is based on a different correction score in comparison to the original debiased Lasso. The specific expression and rationale are illustrated in detail in Section \ref{sec2-2}. In Step (ii), the estimator $\hat{\beta}^*$ is a noiseless Lasso based on the empirical estimates. In the usual parametric bootstrap, the response vector is constructed as
\begin{align}
\label{bs1}
  y^*_i = x_i\hat{\beta}+ \hat{\sig}\xi^*_i,~ i= 1,\dots,n,
\end{align}
where $\xi^*_i$ are synthetic \textit{i.i.d.} standard Gaussian random variables and $\hat{\sig}$ is a consistent estimate of noise level. Hence, the estimator $\hat{\beta}^*$ employed in Step (ii) can be viewed as a noiseless Lasso based on parametric bootstrap. In fact, our $\hat{\beta}^*$ can be replaced with an average of the usual noisy parametric bootstrap estimators. Although the noisy parametric bootstrap estimates can be used to generate empirical confidence intervals, our main purpose of bootstrap is for bias correction rather than uncertainty quantification. Hence, we focus on the proposed noiseless version which can simplify the computation and does not introduce extra randomness. (I would like to thank an anonymous referee for pointing this out.)

We also comment that the current framework and the proposed approach are different from the results of bootstrapping the Adaptive Lasso considered in \citet{CL13}. In fact, in the current work, the beta-min condition is not required and selection consistency is not needed. As a consequence, some direct methods, such as least squares after selection, cannot lead to valid inference procedures under current conditions in general.

\subsection{The proposed correction score}
\label{sec2-2}
Before articulating the format of the proposed correction score $\hat{w}_j$, let us bring some intuitions to its construction. Recall that the success of the debiased Lasso relies on choosing a correction score $v\in\R^n$ such that \citep{ZZ14}
\begin{equation}
\label{eq1-4-1}
e_j^T-v^TX/n\approx 0.
\end{equation}
One possible realization of $v$ is $\hat{z}_j$ defined in (\ref{eq-z1}). In fact, it is ideal but not amendable to have an ``exactly equal'' relationship in (\ref{eq1-4-1}). It can be seen from later discussion (\ref{err1}) that if $v$ is such that  $v^TX_A/n= 0$ for a set $A\subseteq \{-j\}$, then the debiased Lasso based on the correction score $v$ is free from the bias of $\hat{\beta}_{A}$. However, an ``exactly equal'' relationship is not achievable in general since $\Theta^n=(\Sig^n)^{-1}$ is not well-defined when $p\geq n$. As a compromise, we would like to obtain a correction score $v$ such that
\[
   (e_j^T-v^TX/n)_S= 0~~\text{and}~~ (e_j^T-v^TX/n)_{S^c}\approx 0.
\]
That is, on a ``small'' set which has large bias, $S$, we require the exact equality to hold in (\ref{eq1-4-1}) and on $S^c$, we allow for ``approximately equal''. As the true support $S$ is unknown, we replace it by an estimate based on the Lasso, $\widehat{S}=\{j:\hat{\beta}_j\neq 0\}$, which is not necessarily consistent but satisfies certain desirable properties.  
Let $\widehat{A}_j=\widehat{S}\cap \{-j\}$. Let $P_{\widehat{A}_j}^{\perp}=I_{n}-X_{\widehat{A}_j}\left(X_{\widehat{A}_j}^TX_{\widehat{A}_j}\right)^{-1}X^T_{\widehat{A}_j}$ denote the projection operator, where $X_{\widehat{A}_j}^TX_{\widehat{A}_j}/n$ will be shown to be invertible with high probability. Formally, our proposed correction score $\hat{w}_j$ can be defined as
\begin{equation}
\label{eq-z2}
  \hat{w}_{j}=\left\{\begin{array}{ll}
   \frac{nP_{\widehat{A}_j}^{\perp}X\Theta_j}{\langle X_j,P_{\widehat{A}_j}^{\perp}X\Theta_j\rangle } ~~\text{if $\Theta$ is known}\\
  \frac{nX\hat{\kappa}}{\langle X_j,X\hat{\kappa}\rangle }~ ~\text{if $\Theta$ is unknown},
  \end{array}\right.
\end{equation}
where
\begin{equation}
\label{hkappa}
   \hat{\kappa}=\argmin_{\kappa\in\R^{p}}\left\{\frac{1}{2n}\|X\kappa\|_2^2+\lam_j\|\kappa_{\widehat{S}^c\setminus \{j\}}\|_1:~\kappa_j=1\right\}.
\end{equation}
The optimization in (\ref{hkappa}) is different from (\ref{hgam0}) as it does not penalize the coefficients in $\widehat{A}_j$. 
It is easy to see from (\ref{eq-z2}) and the KKT condition of (\ref{hkappa}) that $(e^T_j-\hat{w}_j^TX/n)_{\widehat{S}}=0$ no matter $\Theta$ is known or not. We demonstrate the theoretical advantages of the proposed BS-DB estimator in the next section.


\section{Theoretical properties}
\label{sec3}
In this section, we study the theoretical properties of the proposed BS-DB estimator. We will first show some preliminary results under the irrepresentable condition and then justify the asymptotic validness of the proposed confidence interval (\ref{ci1}).  We will further study the robustness of our proposal to the violation of the irrepresentable condition. For the main results, we assume the following conditions.
\begin{condition}[Gaussian designs]
\label{cond1b}
Each row of $X$ is i.i.d. Gaussian distributed with mean 0 and covariance matrix $\Sig$, where $\Sigma$ satisfies that $0<C_{\min}\leq \Lambda_{\min}(\Sig)\leq \Lambda_{\max}(\Sig)\leq C_{\max}<\infty$ and $\max_{1\leq j\leq p} \Sig_{j,j}\leq K_1<\infty$.
\end{condition}
The smallest eigenvalue of $\Sig$ is assumed to be bounded away from zero, which is used in upper bounding the remaining bias of the debiased Lasso. The largest eigenvalue of $\Sig$ is assumed to be finite, which is used in justify the asymptotic normality. The Gaussian distribution is needed for technical convenience and is also required in closely related works \citep{Wain09, JM15}. 
\begin{condition}[Sub-Gaussian errors]
\label{cond2b}
$\eps_i~(i=1,\dots,n)$ are \textit{i.i.d.} sub-Gaussian with mean 0 and variance $\sig^2>0$ such that $\max_{1\leq i\leq n}\|\eps_i\|_{\psi_2}\leq K_2$.
\end{condition}
\begin{condition}[Irrepresentable condition]
\label{cond3b}
\[
\phi=\max_{k\in S^c}\left\|\Sig_{k,S}(\Sig_{S,S})^{-1}\right\|_{1} < 1.
\]
\end{condition}
As previously mentioned, Condition \ref{cond3b} is crucial in our main analysis. It has been used in analyzing the variable selection properties of the Lasso \citep{ZY06, Wain09} and it is equivalent to the ``neighborhood stability'' condition \citep{MB06}. In fact, the purpose of irrepresentable condition is to guarantee that the support estimated by the Lasso is a subset of the true support with high probability (Lemma \ref{lem1a} and Lemma \ref{lem1b}). Some concave penalized methods, say \citet{MCP}, require weaker versions of Condition \ref{cond3b} for guaranteeing such properties and they can also be suitable for debiasing. 
For conciseness and computational convenience, we focus on the Lasso. 

It is worth mentioning that Condition \ref{cond3b} does not rule out the most difficult scenarios which yield the minimax optimal sample size condition  $n\gg (s\log p)^2$. Indeed, the lower bound results for the confidence intervals in high-dimensional linear models \citep{CG17,JM15} are based on some worst cases satisfying $\phi=0$. Therefore, our results are comparable to the existing lower bound results achieved by the debiased Lasso.


\subsection{Preliminary lemmas}
\label{sec3-1}
We first prove some preliminary results as consequences of the irrepresentable condition.
For $b\in\R^p$, let $sgn(b)$ be an element of the sub-differential of the $\ell_1$-norm of $b$. 
That is, $\left(sgn(b)\right)_j=sgn(b_j)=b_j/|b_j|$ if $b_j\neq 0$ and $\left(sgn(b)\right)_j\in[-1,1]$ if $b_j=0$.  Let 
\begin{equation}
\label{rho-s}
 \rho_S=\max_{j\in S} \left\|e_j^T(\Sig_{S,S})^{-1}\right\|_1.
 \end{equation}
Notice that $\rho_S\in[c_1, c_2\sqrt{s}]$ for some positive constants $c_1$ and $c_2$ under Condition \ref{cond1b}.
 \begin{lemma}[Sparsity pattern recovery for the Lasso]
\label{lem1a}
Assume that Condition \ref{cond1b} - Condition \ref{cond3b} are satisfied.
If the sequence $(n,p,s,\lam)$ satisfies that 
\begin{align}
\label{lem1a-n}
 n\geq \frac{cK_1s\log p}{C_{\min}(1-\phi)^2}~~\text{and}~\lam\geq \frac{c_0\sig}{1-\phi}\sqrt{\frac{K_1\log p}{n}}
 \end{align}
for large enough constants $c$ and $c_0$, then 
\begin{align}
&\Omega_0:=\{\widehat{S}\subseteq S~\text{and}~\|\hat{\beta}_S-\beta_S\|_{\infty}\leq C\rho_S\lam \}\label{omega0}
\end{align}
holds with probability at least $1-c_1/p-\exp(-c_2n)-c_3/n$ for some large enough constants $C$, $c_1$, $c_2$, $c_3$, and $\rho_S$ defined in (\ref{rho-s}).
\end{lemma}
Lemma \ref{lem1a} is more general than Theorem 3 in \citet{Wain09} in the sense that the $\ell_{\infty}$-bound of the estimation error in (\ref{omega0}) does not require the beta-min condition. This is a nontrivial generalization as $sgn(\hat{\beta}_S)\neq sgn(\beta_S)$ in general and dependence between $sgn(\hat{\beta}_S)$ and $X_S$ makes it challenging to get the desired bound on $\hat{\beta}_S-\beta_S$. Our idea is to decouple $\hat{\beta}_{S}$ with $X_k$ for each $k\in S$ with a ``leave-one-out'' argument. 

 The sparsity pattern recovery can be similarly derived for the bootstrap version $\hat{\beta}^*$. We summarize an important consequence of Lemma \ref{lem1a} and its bootstrap counterpart in the next lemma.
 Let $\widehat{S}^*=\{j: \hat{\beta}^*_j\neq 0\}$.
\begin{lemma}
\label{lem1b}
Assume that Condition \ref{cond1b} - Condition \ref{cond3b} are satisfied.
If (\ref{lem1a-n}) holds, then 
\begin{align}
\label{omega0*}
& \Omega_0^*:=\{\widehat{S}^*\subseteq \widehat{S}~\text{and}~\|\hat{\beta}^*_{S}-\hat{\beta}_S\|_{\infty}\leq C^*\rho_S\lam\}.
\end{align}
holds with probability at least $1-c_1/p-c_2/n-\exp(-c_3n)$ for some large enough constants $c_1$, $c_2$, $c_3$, and $C^*$.
 In the event $\Omega_0\cap \Omega_0^*$, $S_+\subseteq \widehat{S}\subseteq S$, $S_+\subseteq \widehat{S}^*\subseteq S$, and $sgn(\beta_j) = sgn(\hat{\beta}_j)= sgn(\hat{\beta}^*_j)$ for $j\in S_+$, where
\begin{equation}
\label{eq-sp}
   S_+=\left\{j:|\beta_j|>  (C+C^*)\rho_S\lam, j=1,\dots,p\right\}
\end{equation}
for $C$ in (\ref{omega0}) and $C^*$ in (\ref{omega0*}).
\end{lemma}
Lemma \ref{lem1b} implies that the estimated supports based on $\hat{\beta}$ and $\hat{\beta}^*$ can be upper and lower bounded by two unknown but deterministic sets, namely $S$ and $S_+$, with high probability, where $S_+$ can be viewed as the set of strong signals. In addition, the sign vectors of $\hat{\beta}_{S_+}$ and $\hat{\beta}^*_{S_+}$ are asymptotically deterministic. These observations are crucial for the success of our proposed bias correction as will be seen from the next subsection. 

\subsection{Illustrations of the proposed bias correction}
\label{sec3-2}
In this subsection, we layout some key steps for analyzing the proposed approach based on the results proved in Section \ref{sec3-1}. Let us focus on the case where $\Theta$ is unknown. We first review some key steps in the typical analysis of the debiased Lasso. From (\ref{beta-db}) and (\ref{omega0}), we have
\begin{align*}
    \hat{\beta}^{(DB)}_j-\beta_j&=\underbrace{\left(e_j^T-\hat{z}_j^TX/n\right)(\hat{\beta}-\beta)}_{\widehat{rem}_j^o}+\frac{\langle \hat{z}_j,\eps\rangle }{n},
    \end{align*}
    where $\widehat{rem}^o_j=O_P(s\log p/n)$ by an $\ell_{\infty}-\ell_1$ splitting using the KKT condition of (\ref{hgam0}) and $\ell_1$-bound on $\hat{\beta}-\beta$. 
For the asymptotic normality of $\sqrt{n}(\hat{\beta}^{(DB)}_j-\beta_j)$ to be true, one needs $|\widehat{rem}^o_j|=o_P(n^{-1/2})$, which gives the typical sample size condition $(s\log p)^2\ll n$, and the asymptotic normality of $\langle \hat{z}_j,\eps\rangle/\sqrt{n}$, which holds under mild conditions.
   
Next, we illustrate that, with the proposed BS-DB, the bias coming from strong signals is removed in two steps. Let $W^n =X^T\eps/n$. The event $\Omega_0$ (\ref{omega0}) and the KKT condition of $\beta$ together imply that
\begin{align}
    \hat{\beta}^{(mDB)}_j-\beta_j&=\left(e_j^T-\hat{w}_j^TX/n\right)_S(\hat{\beta}_S-\beta_S)+\frac{\langle \hat{w}_j,\eps\rangle }{n} \nonumber\\
    &=\underbrace{-\lam\left(e_j^T-\hat{w}_j^TX/n\right)_S (\Sig^{n}_{S,S})^{-1}\sgn(\hat{\beta}_S)}_{\hat{b}_j}\nonumber\\
    &\quad+\underbrace{\left(e_j^T-\hat{w}_j^TX/n\right)_S (\Sig^{n}_{S,S})^{-1}W^n_S}_{\hat{r}_j}+\frac{\langle \hat{w}_j,\eps\rangle }{n},\label{err1}
    \end{align}
    where $\hat{b}_j$ is the dominant bias term and $\hat{r}_j$ and the last term are linear combinations of $\eps$.  For $S_+$ defined in Lemma \ref{lem1b}, let $S_-=S\setminus S_+$ denote the set of weak signals. Let $s_+=|S_+|$ and $s_-=|S_-|$. To see the bias contributed by strong signals, we rewrite $\hat{b}_j$ as
    \begin{align*}
    \hat{b}_j=    -\lam\left(e_j^T-\hat{w}_j^TX/n\right)_S\underbrace{\begin{pmatrix}
(\Sig^{n}_{S,S})_{S_+,S_+}^{-1}\sgn(\hat{\beta}_{S_+}) &(\Sig^{n}_{S,S})_{S_+,S_-}^{-1}\sgn(\hat{\beta}_{S_-})\\
(\Sig^{n}_{S,S})_{S_-,S_+}^{-1}\sgn(\hat{\beta}_{S_+}) & (\Sig^{n}_{S,S})_{S_-,S_-}^{-1}\sgn(\hat{\beta}_{S_-})
    \end{pmatrix}}_{G^n},
    \end{align*}
    where $\sgn(\hat{\beta}_{S_+})=\sgn(\beta_{S_+})$ with high probability under the conditions of Lemma \ref{lem1b}. We will denote $G^n$ as a $2\times 2$ block matrix such that $G^n_{1,1}=(G^n)_{S_+,S_+}$.
   Loosely speaking, the purpose of bootstrap is to remove $G^n_{1,1}$ and $G^n_{2,1}$ and the purpose of proposed correction score $\hat{w}_j$ is to remove $G^n_{2,1}$. Specifically, in event $\Omega^*_0$ (\ref{omega0*}), we can similarly show that for $ \hat{b}^*_j$ defined in (\ref{eq-hbj}),
    \begin{align*}
 \hat{b}^*_j&=-\lam\left(e_j^T-\hat{w}_j^TX/n\right)_S\underbrace{\begin{pmatrix}
(\Sig^{n}_{S,S})_{S_+,S_+}^{-1}\sgn(\hat{\beta}^*_{S_+}) &(\Sig^{n}_{S,S})_{S_+,S_-}^{-1}\sgn(\hat{\beta}^*_{S_-})\\
(\Sig^{n}_{S,S})_{S_-,S_+}^{-1}\sgn(\hat{\beta}^*_{S_+}) & (\Sig^{n}_{S,S})_{S_-,S_-}^{-1}\sgn(\hat{\beta}^*_{S_-})
    \end{pmatrix}}_{G^*}.
\end{align*}

    Since $sgn(\hat{\beta}^*_{S_+})=sgn(\hat{\beta}_{S_+})$ by Lemma \ref{lem1b}, we have $G^n_{.,1}=G^*_{.,1}$. Moreover, the second column of $G^n$ and that of $G^*$ come from the incorrect sign estimation of weak signals. The term $G^n_{1,2}$ is caused by the possible correlation between $X_{S_+}$ and $X_{S_-}$ and it can be large if $s_+$ is large. 
 To get rid of $G^n_{1,2}$ and $G^*_{1,2}$, we invoke that the proposed correction score $\hat{w}_j$ satisfies $(e_j^T-\hat{w}_j^TX/n)_{\widehat{S}}=0$. If $S_+\subseteq \widehat{S}$, then the effects of $G^n_{1,2}$ and $G^*_{1,2}$ are removed by the proposed correction score. To summarize, (\ref{err1}) and above analysis together imply that in $\Omega_0\cap\Omega_0^*$, 
   \begin{align*}
    &\hat{\beta}^{(mDB)}_j-\beta_j-\hat{b}^*_j= -\lam\left(e_j^T-\hat{w}_j^TX/n\right)_{S_-}(G_{2,2}^n  -G^*_{2,2})  +\hat{r}_j+\frac{\langle \hat{w}_j,\eps\rangle }{n}.
        \end{align*}
We will show that the remaining bias, i.e. the first two terms on the right hand side on the above expression, are $O_P(s_-\log p/n)$. 
    Finally, the validness of the confidence interval considered in (\ref{ci1}) also relies on the asymptotic normality of $\langle \hat{w}_j,\eps\rangle /\sqrt{n}$. We will employ the central limit theorem to show desirable results under proper conditions.

\subsection{Main theorems}
\label{sec3-3}
In this subsection, we formally establish the main theorems for statistical inference with the BS-DB approach. 
As a benchmark, we first present an error analysis of the original debiased Lasso demonstrating its different magnitude of bias on and off the true support.
\begin{theorem}[The remaining bias of the debiased Lasso]
\label{thm2a}
Suppose that Condition \ref{cond1b} - Condition \ref{cond3b} hold true and $s\log p\ll n$. 
Let $\hat{\beta}^{(DB)}_j$ be computed via (\ref{beta-db}) with $\lam=c_1\sqrt{\log p/n}$ and $\lam_j=c_2\sqrt{\log p/n}$ for some sufficiently large constants $c_1$ and $c_2$ such that (\ref{lem1a-n}) holds. Then
 \[
   \hat{\beta}_j^{(DB)}-\beta_j=\frac{\langle \hat{z}_j,\eps\rangle}{n}+ \widehat{rem}^o_j,
\]
where 
\begin{align}
\label{eq-re1}
|\widehat{rem}^o_j|=_{O_p}\left\{
\begin{array}{ll}
 \frac{\sqrt{s\log p}}{n}+ \frac{\rho_S \sqrt{s\log p}}{n}\mathbbm{1}(j\in S)~~\text{if $\Theta$ is known}\\
 \frac{s\log p}{n}~~\text{if}~\Theta~\text{is unknown}.
\end{array}\right.
\end{align}
\end{theorem}
Theorem \ref{thm2a} implies that when $\Theta$ is known, the remaining bias of $\hat{\beta}_j^{(DB)}$ can be  larger for $j\in S$ than for $j\in S^c$ as $1\lesssim \rho_S\lesssim \sqrt{s}$ under Condition \ref{cond1b}. Moreover, as $\rho_S$ gets larger, the bias of $\hat{\beta}_j^{(DB)}$ for $j\in S$ gets larger.
 These results coincide with our observations in Figure \ref{fig0}.
We mention that one can get rid of the term which involves $\rho_S$ by estimating $\hat{\beta}$ and $\hat{z}_j$ with two independent subsets of data (Theorem 7 in \citet{CG17} and Proposition H.1. in \citet{JM15}). In practice, especially when the sample size is relatively small, sample splitting can produce unstable results. 
Hence, we focus on the version without sample splitting.
When $\Theta$ is unknown, the results of Theorem \ref{thm2a} agree with  the existing analysis of the debiased Lasso \citep{ZZ14, Van14, JM14}. With Theorem \ref{thm2a} serving as a benchmark, we study the  the remaining bias of the proposed BS-DB in the next lemma.

 \begin{lemma}[The remaining bias of BS-DB]
  \label{lem2b}
 Suppose that Condition \ref{cond1b} - Condition \ref{cond3b} hold true and $s\log p\ll n$. Let $\hat{\beta}^{(BS-DB)}_j$ be computed via (\ref{beta-bsdb}) with $\lam=c_1\sqrt{\log p/n}$ and $\lam_j=c_2\sqrt{\log p/n}$ for some sufficiently large constants $c_1$ and $c_2$ such that (\ref{lem1a-n}) holds. Then the expressions in (\ref{eq-bsdb0}) and (\ref{eq-re2}) hold true.
 \end{lemma}
 Lemma \ref{lem2b} shows that the magnitude of the remaining bias of $\hat{\beta}^{(BS-DB)}_j$ is determined by the number of weak signals. Comparing with Theorem \ref{thm2a}, we see that the remaining bias of $\hat{\beta}^{(BS-DB)}_j$ is much smaller than the remaining bias of the debiased Lasso when $s_-\ll s$. This  demonstrates the improvement of our proposal and convinces our observations in Figure \ref{fig01}. When $\Theta$ is known, the magnitude of remainder terms of $\hat{\beta}_j^{(BS-DB)}$ can be different for $j\in S$ and $j\notin S$, which is analogous to the results for debiased Lasso.  

Next, we move on to establish the limiting distribution of the BS-DB estimator. We first prove the convergence rate of a variance estimator, which can also benefit from the irrepresentable condition. Define
\begin{align}
\label{sig-est}
\hat{\sig}^2= \frac{1}{n-|\widehat{S}|}\|P_{\widehat{S}}^{\perp}y\|_2^2,
\end{align}
where $P_{\widehat{S}}^{\perp}=I_{n}-X_{\widehat{S}}\left(X_{\widehat{S}}^TX_{\widehat{S}}\right)^{-1}X^T_{\widehat{S}}$.

\begin{lemma}[Convergence rate of the variance estimator]
 \label{lem0b}
 Suppose that Condition \ref{cond1b} - Condition \ref{cond3b} hold true and $s\log p\ll n$. Let $\hat{\beta}$ be computed via (\ref{beta-init}) with $\lam=c_1\sqrt{\log p/n}$ for a sufficiently large constant $c_1$ satisfying (\ref{lem1a-n}). For $\hat{\sig}^2$ defined in (\ref{sig-est}), 
 \begin{align}
  \left|\hat{\sig}^2-\sig^2\right|=O_P\left(n^{-1/2}+\frac{\min\{\rho_S^2s_-,s\}\log p}{n}\right).\label{lem0b-re1}
 \end{align}
 \end{lemma}
 On the right hand side of (\ref{lem0b-re1}), the first term $n^{-1/2}$ comes from the randomness of $\eps$ and the second term comes from the estimation error of weak signals. It is known that a widely used variance estimator, the mean of squared residuals based on $\hat{\beta}$, has convergence rate $O_P(n^{-1/2}+s\log p/n)$. We see that $\hat{\sig}^2$ is no worse than and can have faster rate of convergence than the mean of squared residuals based on the Lasso if $\rho_S^2s_-\ll s$. That is, $\hat{\sig}^2$ can be more accurate when the number of strong signals is dominant and the correlation among $X_S$ is weak. The empirical performance of $\hat{\sig}^2$ is evaluated in various numerical experiments in Section \ref{sec-simu}.

In the next theorem, we collect all the preliminary results and prove the asymptotic normality of $\hat{\beta}^{(BS-DB)}_j$ under proper sample size conditions. Let $\Phi(c)=\P(x\leq c)$ for a standard normal variable $x$.
\begin{theorem}[Asymptotic normality of BS-DB]
\label{thm2b}
Suppose that Condition \ref{cond1b} - Condition \ref{cond3b} hold true and $s\log p\ll n$. Let $\hat{\beta}^{(BS-DB)}_j$ be computed via (\ref{beta-bsdb}) with $\lam=c_1\sqrt{\log p/n}$ and $\lam_j=c_2\sqrt{\log p/n}$ for some sufficiently large constants $c_1$ and $c_2$ such that (\ref{lem1a-n}) holds. (i) If $\Theta_{j}$ is known, then
\begin{align}
\label{thm2b-eq1}
\sup_{c\in \R}\left|\P\left( \hat{\beta}_j^{(BS-DB)}-\beta_j\leq c\hat{\sig}\|\hat{w}_j\|_2/n\right)- \Phi(c)\right| \rightarrow 0.
\end{align}
(ii) When $\Theta_{j}$ is unknown, (\ref{thm2b-eq1}) holds if $\max\left\{s_-,\|\Theta_{(S_+)^c,j}\|_0\right\}\log p\ll \sqrt{n}$.
\end{theorem}

Theorem \ref{thm2b} implies that the proposed confidence interval (\ref{ci1}) has nominal coverage probability asymptotically. The condition $s_-\log p\ll \sqrt{n}$ when $\Theta$ is unknown guarantees that the remaining bias of BS-DB is $o_P(n^{-1/2})$. The condition $\|\Theta_{(S_+)^c,j}\|_0\log p\ll\sqrt{n}$ guarantees that $\hat{\kappa}$ converges to its probabilistic limit at a sufficiently fast rate such that $\hat{\eta}_j$ in (\ref{eq-bsdb0}) is asymptotically normal. It is especially needed here as our constructed $\hat{w}_j$ is dependent with $\eps$. This condition can be relaxed if $\widehat{S}$ used in (\ref{eq-z2}) is independent of $X$ or is asymptotically deterministic. Hence, one can perform a sample splitting and compute $\widehat{S}$ with one fold of the data and compute $\hat{\kappa}$ and other estimates with the other fold of the data. In this way, the sparsity requirement on $\Theta_{(S_+)^c,j}$ can be relaxed. On the other hand, some mild conditions can lead to asymptotically deterministic $\widehat{S}$. One sufficient condition is that $|\{j:~|\beta_j|\asymp \rho_S\lam\}|=0$ which guarantees $\widehat{S}=S_+$ asymptotically.
The \textbf{optimality} of the proposed confidence interval (\ref{ci1}) can be partially understood from the established lower bound results. Specifically, the parameter space considered in Theorem 3 of \citet{CG17} is
 \[
   \Xi_0(s)=\{(\beta,\sig,\Sig):\|\beta\|_0\leq s,\sig^2\geq c_1, 1/c_2\leq \Lambda_{\min}(\Sig)\leq \Lambda_{\max}(\Sig)\leq c_2 \}
 \]
 for some constants $c_1>0$ and $c_2>1$. 
  Let $L^*_{\alpha}(\Theta,\beta_j)$ denote the minimax expected length of confidence intervals for $\beta_j$ ((2.3) of \citet{CG17}) over $\Theta$ at confidence level $\alpha$. Recall that $\underline{s}=|\{j:|\beta_j|\asymp \sqrt{\log p/n}\}|$. When $\Theta$ is unknown, it has been shown that
 \[
    L^*_{\alpha}(\Xi_0(s),\beta_j)\gtrsim \frac{1}{\sqrt{n}}+\frac{\underline{s}\log p}{n}.
 \]
 Let us consider a more detailed parameter space
 \begin{align*}
     \Xi(s_+,s_-)=&\{(\beta,\sig,\Sig):|S_+|\leq s_+,|S_-|\leq s_-,\sig^2\geq c_1, \\
    & ~~1/c_2\leq \Lambda_{\min}(\Sig)\leq \Lambda_{\max}(\Sig)\leq c_2,~\rho_S\leq c_3\},
 \end{align*}
 for some constants $c_1>1$, $c_2>1$, and $c_3<\infty$. As $s_-=\underline{s}$ for constant $\rho_S$, it is not hard to see that
 \[
    L^*_{\alpha}(\Xi(s_+,s_-),\beta_j)\gtrsim \frac{1}{\sqrt{n}}+\frac{s_-\log p}{n}.
 \]
This shows the optimality of the proposed BS-DB in $\Xi(s_+,s_-)$ when $\Theta$ is unknown.

\subsection{When irrepresentable condition does not hold}
\label{sec3-5}
The irrepresentable condition is always hard to check in reality. Therefore, it is important to understand whether the proposed BS-DB is still valid for inference when such a condition is not true. In the next theorem, we justify the theoretical properties of BS-DB without the irrepresentable condition. We also relax the Gaussian assumption in Condition \ref{cond1b} to sub-Gaussian designs. For practical concerns, we only prove for the case where $\Theta$ is unknown.
\begin{condition}[Sub-Gaussian designs]
\label{cond4b}
Each row of $X$ is i.i.d from a sub-Gaussian distribution with mean zero and covariance matrix $\Sig$. The matrix $\Sigma$ satisfies that $0<C_{\min}\leq \Lambda_{\min}(\Sig)\leq \Lambda_{\max}(\Sig)\leq C_{\max}<\infty$. There exists a positive constant $K_1$ such that $\max_{1\leq i\leq n} \|x_i\Sig^{-1/2}\|_{\psi_2}\leq K_1<\infty.$
\end{condition}
\begin{theorem}[The remaining bias of BS-DB without irrepresentable condition]
\label{thm10}
Assume that Condition \ref{cond2b} and Condition \ref{cond4b} are true and $s\log p\ll n$. Let $\hat{\beta}^{(BS-DB)}_j$ be computed via (\ref{beta-bsdb}) with $\lam=c_1\sqrt{\log p/n}$ and $\lam_j=c_2\sqrt{\log p/n}$ for some sufficiently large constants $c_1$ and $c_2$. When $\Theta$ is unknown, it holds that
\[
  \hat{\beta}^{(BS-DB)}_j-\beta_j=\frac{\langle \hat{w}_j,\eps\rangle}{n}+O_P\left(\frac{s\log p}{n}\right).
\]
\end{theorem}
Theorem \ref{thm10} implies that if $n\gg (s\log p)^2$, then the remaining bias of $\hat{\beta}_j^{(BS-DB)}$ is $o_P(n^{-1/2})$ even if the irrepresentable condition does not hold. Hence, the remaining bias of BS-DB estimator is no larger than that of the debiased Lasso. The asymptotic normality of $\langle \hat{w}_j,\eps\rangle/\sqrt{n}$ can be established based on the proof of Theorem \ref{thm2b} under proper conditions.
We can summarize that when $\Theta$ is unknown, there is no loss in applying the BS-DB asymptotically regardless of the irrepresentable condition and BS-DB can achieve more accurate confidence intervals in existence of a large proportion of strong signals.
\section{Numerical experiments}
\label{sec-simu}
In this section, we demonstrate the empirical performance of BS-DB in comparison to the debiased Lasso in more practical settings with $\Sig$ and $\sig^2$ are both unknown. 

We set sample size $n=100$, the number of covariates $p=300$, and the noise level $\sig=1$. We consider $\Sig=I_p$, $\Sig=\Sig^o$ defined in Section \ref{sec1-3}, and another $\Sig$ where the irrepresentable condition does not hold. Each row of $X$ is \textit{i.i.d.} generated from $N(0,\Sig)$ and $\eps_i\sim N(0,\sig^2)$ for $i=1,\dots,n$. We consider three levels of sparsity. Specifically, $s=4k$ for $k\in\{1,2,3\}$. Let $s_-=|\{j: \beta_j=0.2\}|$. For any $k\in\{1,2,3\}$, we consider the following two cases
\begin{enumerate}
\item[(i)] $(s,s_-)=(4k,1)$, $\beta_{1:s}=(4,2,4,2,\dots, 4,0.2)^T$, and $\beta_{s+1:p}=\mathbf{0}_{p-s}$.
\item[(ii)] $(s,s_-)=(4k,2k)$, $\beta_{1:2k}=(4,2,\dots, 4,2)^T$, $\beta_{2k+1:s}=(0.2,\dots,0.2)^T$, and $\beta_{s+1:p}=\mathbf{0}_{p-s}$.
\end{enumerate}
In case (i), a large proportion of signals are strong and in case (ii), half of the signals are strong.
The purpose is to demonstrate the effect of overall sparsity and effect of the proportion of weak signals separately. We construct two-sided 95\% confidence intervals for $\beta_1=4$, $\beta_2=2$, $\beta_s=0.2$, and $\beta_{s+1}=0$ in each setting. We report summarized statistics based on 500 independent realizations for each setting. 

We will compare the coverage probabilities and lengths of confidence intervals  given by the debiased Lasso (DB) procedure and the proposed BS-DB. 
In both methods, the tuning parameter for the Lasso is $\lam=2\hat{\sig}^{(init)}\sqrt{\log p/n}$, where $\hat{\sig}^{(init)}$ is computed via the scaled Lasso. The tuning parameter in (\ref{hgam0}) and (\ref{hkappa}) are set to be $\hat{\sig}_x^{(init)}\sqrt{2\log p/n}$, where $\hat{\sig}_x^{(init)}$ is computed via the scaled Lasso with response $X_j$ and covariates $X_{-j}$. The estimated noise level $\hat{\sig}$ is computed according to (\ref{sig-est}).

\subsection{Identify covariance matrix}
\label{sec-simu1}
We first present the numerical results for the case where $\Sig=I_p$ (Table \ref{table1}). 
One can see that the coverage probabilities of debiased Lasso are less accurate than those of BS-DB in all the scenarios, especially when the number of strong signals dominant. In comparison, BS-DB has coverage probabilities close to the nominal level at different sparsity levels and is especially robust to the existence of a large proportion of strong signals. This agrees with our analysis in Theorem \ref{thm2a} and Theorem \ref{thm2b}. The lengths of confidence intervals produced by the debiased Lasso and BS-DB are comparable. 
Hence, the difference in coverage probabilities is mainly due to the remaining bias term. 

The confidence intervals based on debiased Lasso have relatively poor coverage in our numerical studies comparing with in some previous studies. 
The main reason is that previous studies always use mean of squared residuals (MSR) to estimate the noise level which tends to inflate the true value and results in wider confidence intervals, which compensate for the bias of debiased Lasso.
The performance of the proposed variance estimator (\ref{sig-est}), the scaled Lasso estimator \citep{scaled-lasso}, and the MSR are reported in Table \ref{table-sig}. One can see that the proposed estimator of noise level is most accurate in different settings. The MSR can have large bias when the number of strong signals dominants. The scaled-Lasso is more reliable than MSR but not as accurate as the proposed estimator when the sparsity is large.

\begin{table}[ht]
\centering
\begin{tabular}{cc|cccc|cccc}
\hline
  & & \multicolumn{4}{c|}{$s_-=1$} &  \multicolumn{4}{c}{$s_-=s/2$} \\
  \hline
  $s$ & $\beta$ & cov.bsdb  & cov.db  & se.bsdb & se.db & cov.bsdb  & cov.db & se.bsdb & se.db \\ 
  \hline
 4 & 4 & 0.942  & 0.888 & 0.10 & 0.10 & 0.936 & 0.890  & 0.10 & 0.10 \\ 
   4 & 2 & 0.938  & 0.882 & 0.10 & 0.10 & 0.932 & 0.902  & 0.10 & 0.10 \\ 
   4 & 0.2 & 0.918  & 0.866 & 0.10 & 0.10 & 0.948  & 0.926 & 0.10 & 0.10 \\ 
   4 & 0 & 0.958 & 0.828 & 0.10  & 0.10 & 0.948 & 0.886  & 0.10& 0.10 \\ 
   \hline
  8 & 4 & 0.932 & 0.666   & 0.11& 0.11 & 0.900  & 0.838 & 0.11 & 0.11 \\ 
   8 & 2 & 0.950 & 0.652  & 0.11 & 0.11 & 0.934 & 0.800  & 0.11 & 0.11 \\ 
 8 & 0.2 & 0.940  & 0.696 & 0.11 & 0.11 & 0.962 & 0.840  & 0.11 & 0.11 \\ 
  8 & 0 & 0.942  & 0.666 & 0.11 & 0.11 & 0.942  & 0.818 & 0.11 & 0.11 \\ 
  \hline
 12 & 4 & 0.936 & 0.498 & 0.11 & 0.11 & 0.946  & 0.742 & 0.11 & 0.11 \\ 
 12 & 2 & 0.922 & 0.488  & 0.11 & 0.11 & 0.932  & 0.716 & 0.11 & 0.11 \\ 
   12 & 0.2 & 0.948  & 0.544 & 0.11 & 0.11 & 0.956  & 0.764 & 0.11 & 0.11 \\ 
  12 & 0 & 0.956  & 0.552 & 0.11 & 0.11 & 0.962 & 0.788 & 0.11  & 0.11 \\ 
   \hline
\end{tabular}
\caption{Coverage probabilities of BS-DB (cov.bsdb), coverage probabilities of debiased Lasso (cov.db), standard errors of BS-DB (se.bsdb), and standard errors of debiased Lasso (se.db) when $\Sig=I_p$.}
\label{table1}
\end{table}

%
%
\subsection{Mild correlation on the support}
\label{sec4-2}
In this subsection, we consider $\Sig=\Sig^o$ specified in Section \ref{sec1-3}. This is a harder scenario than $\Sig=I_p$ since $\rho_S$ is increased. Other parameters are set to be the same as in Section \ref{sec-simu1}. We see from Table \ref{table2} that the debiased Lasso has coverage probabilities much lower than the nominal level at all the sparsity levels while the BS-DB remains to be robust at most sparsity levels. The most difficult case is $(s,s_-)=(12,1)$, i.e. the overall sparsity is large and most of them are strong signals, where the confidence intervals given by BS-DB provide lower coverage probabilities for strong signals. One reason is that the estimated noise level given by the scaled-Lasso has large errors in this case. Hence, the choice of $\lam$ can be improper and cause large finite sample bias for the initial Lasso estimator. One way to get around this issue is to select $\lam$ via cross validation.
We can see from Table \ref{table-sig} that the proposed variance estimator has the most reliable performance  when $\Sig=\Sig^o$ in comparison to other two estimators.

\begin{table}[ht]
\centering
\begin{tabular}{cc|cccc|cccc}
\hline
  & & \multicolumn{4}{c|}{$s_-=1$} &  \multicolumn{4}{c}{$s_-=s/2$} \\
  \hline
  $s$ & $\beta$ & cov.bsdb & cov.db & se.bsdb & se.db & cov.bsdb &  cov.db & se.bsdb &se.db \\ 
  \hline
 4 & 4 & 0.936 & 0.762 & 0.10 & 0.10 & 0.948  & 0.818  & 0.10 & 0.10 \\ 
 4 & 2 & 0.950 & 0.674& 0.10  & 0.10 & 0.928 & 0.808  & 0.10 & 0.10 \\ 
 4 & 0.2 & 0.944 & 0.662 & 0.10 & 0.10 & 0.924  & 0.786 & 0.10 & 0.10 \\ 
   4 & 0 & 0.946  & 0.846 & 0.10 & 0.10 & 0.940  & 0.876  & 0.10& 0.10 \\ 
   \hline
   8 & 4 & 0.928  & 0.530 & 0.11 & 0.10 & 0.926 & 0.688  & 0.11 & 0.11 \\ 
   8 & 2 & 0.930  & 0.464 & 0.11 & 0.10 & 0.946 & 0.706  & 0.11 & 0.11 \\ 
   8 & 0.2 & 0.958  & 0.492 & 0.11 & 0.10 & 0.926 & 0.702 & 0.11  & 0.11 \\ 
  8 & 0 & 0.952 & 0.646 & 0.11 & 0.10 & 0.948  & 0.810 & 0.11 & 0.11\\ 
  \hline
 12 & 4 & 0.840  & 0.512 & 0.23 & 0.21 & 0.924  & 0.622 & 0.11 & 0.11\\ 
   12 & 2 & 0.904  & 0.378 & 0.23 & 0.21 & 0.936 & 0.610 & 0.11  & 0.11\\ 
   12 & 0.2 & 0.948  & 0.514  & 0.23 & 0.21 & 0.956  & 0.658  & 0.11& 0.11\\ 
  12 & 0 & 0.940 & 0.638 & 0.23 & 0.21 & 0.950 & 0.704  & 0.11 & 0.11\\ 
   \hline
\end{tabular}
\caption{Coverage probabilities of BS-DB (cov.bsdb),  coverage probabilities of debiased Lasso (cov.db), standard errors of BS-DB (se.bsdb), and standard errors of debiased Lasso (se.db) when $\Sig=\Sig^o$.}
\label{table2}
\end{table}

%

\subsection{When the irrepresentable condition does not hold}
In this subsection, we consider a scenario where the irrepresentable condition does not hold. Specifically, we set $\Sig_{j,j}=1$, $1\leq j\leq p$. The upper diagonal of $\Sig$ is such that $\Sig_{j,k}=0.15$ if $s+1\leq k\leq 16$ and $j\leq 16$ and $\Sig_{j,k}=0$ otherwise. The true $\Sig$ remains to be unknown. It is easy to check that for $\phi$ defined in Condition \ref{cond3b}, it holds that $\phi=0.6$ for $s=4$, $\phi=1.2$ for $s=8$ and $\phi=1.8$ for $s=12$. Hence, the irrepresentable condition does not hold with the current $\Sig$ when $s\in\{8,12\}$.
We examine the performance of the BS-DB and debiased Lasso in this case. Other parameters are set to be the same as in Section \ref{sec-simu1}.

In the current setting, the results in Table \ref{table3} show that the BS-DB still provides more accurate coverage than the debiased Lasso but both methods are severely biased when $(s,s_-)=(12,1)$. One reason is that, as in the case of Section \ref{sec4-2}, the estimated noise levels are largely biased and hence the choices of tuning parameters are improper. The confidence intervals are also wider when $(s,s_-)=(12,1)$ and one can see that all three variance estimators are severely biased in this case.
\begin{table}[ht]
\centering
\begin{tabular}{cc|cccc|cccc}
\hline
  & & \multicolumn{4}{c|}{$s_-=1$} &  \multicolumn{4}{c}{$s_-=s/2$} \\
  \hline
  $s$ & $\beta$ & cov.bsdb & cov.db & se.bsdb & se.db & cov.bsdb & cov.db & se.bsdb & se.db \\ 
  \hline
 4 & 4 & 0.934 & 0.886 & 0.10 & 0.10 & 0.916 & 0.908 & 0.10 & 0.10 \\ 
   4 & 2 & 0.908 & 0.850 & 0.10 & 0.10 & 0.932 & 0.910 & 0.10 & 0.10 \\ 
   4 & 0.2 & 0.938 & 0.874 & 0.10 & 0.10 & 0.950 & 0.906 & 0.10 & 0.10 \\ 
   4 & 0 & 0.922 & 0.356 & 0.11 & 0.10 & 0.940 & 0.496 & 0.11 & 0.10 \\ 
   \hline
   8 & 4 & 0.874 & 0.396 & 0.12 & 0.10 & 0.850 & 0.808 & 0.11 & 0.10\\ 
   8 & 2 & 0.874 & 0.390 & 0.12 & 0.10 & 0.856 & 0.758 & 0.11 & 0.10\\ 
   8 & 0.2 & 0.944 & 0.538 & 0.12 & 0.10 & 0.910 & 0.820 & 0.11 & 0.10\\ 
  8 & 0 & 0.926 & 0.006 & 0.12 & 0.10 & 0.826 & 0.076 & 0.11 & 0.10\\ 
  \hline
  12 & 4 & 0.270 & 0.202 & 0.32 & 0.29& 0.820 & 0.620 & 0.11 & 0.10 \\ 
   12 & 2 & 0.326 & 0.162 & 0.33 & 0.29 & 0.814 & 0.608 & 0.11 & 0.10\\ 
   12 & 0.2 & 0.278 & 0.272 & 0.33 & 0.29 & 0.848 & 0.676 & 0.11 & 0.10\\ 
   12 & 0 & 0.084 & 0.000 & 0.33 & 0.29 & 0.630 & 0.002 & 0.11 & 0.10\\ 
   \hline
\end{tabular}
\caption{Coverage probabilities of BS-DB (cov.bsdb),  coverage probabilities of debiased Lasso (cov.db), standard errors of BS-DB (se.bsdb), and standard errors of debiased Lasso (se.db) when the irrepresentable condition does not hold.}
\label{table3}
\end{table}

\begin{table}[ht]
\centering
\begin{tabular}{c|ccc|ccc|ccc}
  \hline
   & \multicolumn{3}{c|}{$\Sig=I_p$} & \multicolumn{3}{c|}{$\Sig=\Sig^o$} &  \multicolumn{3}{c}{IRP does not hold}  \\
  \hline
 $(s,s-)$ & sc-Las & MSR & Prop & sc-Las & MSR & Prop  & sc-Las & MSR & Prop \\ 
  \hline
(4,1) &0.06 & 0.33 & 0.05 &   0.07 & 0.38 & 0.05 & 0.06 & 0.32 & 0.05 \\ 
(4,2) & 0.06 & 0.24 & 0.06 & 0.06 & 0.27 & 0.06 & 0.05 & 0.24 & 0.05\\ 
(8,1) &  0.23 & 0.87 & 0.05 & 0.35 & 0.16 & 0.07   & 0.22 & 0.88 & 0.05\\ 
 (8,4) &  0.12 & 0.50 & 0.07 & 0.14 & 0.56 & 0.06 & 0.09 & 0.47 & 0.05 \\ 
 (12,1) &  0.69 & 1.96 & 0.06 & 1.42 & 3.44 & 0.10  & 1.99& 3.66& 1.90\\ 
 (12,6) & 0.26 & 0.86 & 0.10 & 0.31 & 0.98 & 0.09 & 0.20 & 0.76 & 0.06\\ 
   \hline
\end{tabular}
\caption{Median of absolute errors for estimated noise level with identity covariance matrix, $\Sig=\Sig^o$, and when the irrepresentable condition (IRP) does not hold.}
\label{table-sig}
\end{table}

\section{Discussion}
\label{sec-diss}
In this paper, we propose the BS-DB procedure for constructing confidence intervals for regression coefficients in high-dimensional linear models. Our analysis shows that, under the irrepresentable condition and other mild technical conditions, the BS-DB estimator has smaller order of bias in existence of a large proportion of strong signals in comparison to the debiased Lasso considered in the existing literature. If the irrepresentable condition does not hold, then BS-DB is guaranteed to perform no worse than the debiased Lasso asymptotically. Hence, the BS-DB is a robust and competitive alternative of the original debiased Lasso estimator. 
From our numerical studies, we see that an accurate estimate of noise level is a key ingredient for inference as it can affect the choice of tuning parameters and the uncertainty quantification. It would also be of interest to develop some robust confidence intervals to compensate the bias of the debiased estimators.

\section{Proofs of main lemmas and theorems}
\label{sec-pf}
We first declare some notations. 
 Let $\hat{u}=\hat{\beta}-\beta$ and $\hat{u}^*=\hat{\beta}^*-\hat{\beta}$.
For any $k\in S$, let
\[
   X^{\perp}_k=X_k-\Sig_{k,A_k}\Sig^{-1}_{A_k,A_k}X_{A_k}.
\]
Let $A_k=S\setminus\{k\}$ if $j\in S$ and $A_k=S$ if $j\notin S$. Let $\hat{\beta}^{(k)}_{A_k}$ denote the ``leave-one-out'' Lasso estimate of $\beta_{A_k}$ given the oracle $S$:
\begin{equation}
\label{beta-loo}
  \hat{\beta}^{(k)}=\argmin_{b\in\R^{p}} \left\{ \frac{1}{2n}\|y -X_Sb_S\|_2^2  + \lam\|b_{A_j}\|_1: b_{k}=\beta_k,b_{S^c}=0\right\}.
\end{equation}

We mention that throughout our proof (except for Theorem \ref{thm10}), we repeatedly using the fact that $X\Theta_j$ is independent of $X_{-j}$ when $X$ is Gaussian distributed. This is because the covariance between $X_{-j}$ and $X\Theta_j$ is zero.

\subsection{Proof of Theorem \ref{thm2a}}
\label{pf-lem2a}
\begin{proof}
First consider the case with \underline{\textbf{known $\Theta_j$}}. 
(i) If $j\notin S$,
\begin{align*}
\hat{\beta}^{(DB)}_j-\beta_j=&(e_j^T-\Theta^T_{j}\Sig^n)\hat{u}+\frac{\langle \hat{z}_j,\eps\rangle}{n}.
\end{align*}
In the event that $\{\widehat{S}\subseteq S\}$, $\hat{u}_{S^c}=0$ and (\ref{eq3-2}) holds. As a result,
\begin{align*}
\hat{\beta}^{(DB)}_j-\beta_j&=\underbrace{(-\Theta^T_{j}\Sig^n)_S(\Sig^{n}_{S,S})^{-1}W^{n}_S}_{r_{j,1}} +\underbrace{\lam (\Theta_{j}^T\Sig^n)_S(\Sig^{n}_{S,S})^{-1}sgn(\hat{\beta}_S)}_{r_{j,2}}+\frac{\langle \hat{z}_j,\eps\rangle}{n}.
\end{align*}
Conditioning on $X$, using the sub-Gaussian property of $\eps$, we have with probability at least $1-\exp(-c_1n)$,
\begin{align*}
|r_{j,1}|&\leq C_1K_2\sqrt{\frac{\Theta^T_{j}\Sig^n_{.,S}(\Sig^{n}_{S,S})^{-1}\Sig^n_{S,.}\Theta_{j}}{n}}\leq C_1K_2 \|\Sig^n_{S,.}\Theta_{j}\|_2\|(\Sig^{n}_{S,S})^{-1}\|^{1/2}_2/\sqrt{n}\\
&\leq  C_1K_2\sqrt{s}\|(\Sig^{n}_{S,S})^{-1}\|^{1/2}_2/n,
\end{align*}
where the last step is due to $X\Theta_j$ is independent of $X_S$.
Using Lemma \ref{alem2}, given that $n\gg s$,
\[
   |r_{j,1}|=O_P\left(\frac{K_2\sqrt{K_1s}}{n\sqrt{C_{\min}}}\right).
\]
We bound $r_{j,2}$ as
\begin{align*}
|r_{j,2}|&\leq \lam\left |\Theta_j^T\Sig^n_{.,S}(\Sigma^n_{S,S})^{-1}sgn(\hat{\beta}_{S})\right|\\
&= \lam \left|\Theta_{j,S}sgn(\hat{\beta}_S)+\Theta_{j,S^c}\Sig^n_{S^c,S}(\Sigma^n_{S,S})^{-1}sgn(\hat{\beta}_{S})\right|.
\end{align*}
Define
\begin{equation}
\label{thm1a-eq1}
   X_{S^c}^{\perp}=X_{S^c}-X_{S}\Sig_{S,S}^{-1}\Sig_{S,S^c}.
\end{equation}
Then
\begin{align*}
   (\Sig^{n}_{S,S})^{-1}\Sig^n_{S,S^c}\Theta_{S^c,j}&=\Sig^{-1}_{S,S}\Sig_{S,S^c}\Theta_{S^c,j}+ \frac{1}{n}(\Sig^{n}_{S,S})^{-1}X_S^TX_{S^c}^{\perp}\Theta_{S^c,j}\\
   &=-\Theta_{S,j}+\frac{1}{n}(\Sig^{n}_{S,S})^{-1}X_S^TX_{S^c}^{\perp}\Theta_{S^c,j},
\end{align*}
where the last step is due to $\Theta_{S,j}+ \Sig^{-1}_{S,S}\Sig_{S,S^c}\Theta_{S^c,j}=0$ for $j\notin S$.
Then we have
\begin{align*}
|r_{j,2}|=\lam \left|\frac{1}{n}sgn(\hat{\beta}_S)^T(\Sig^{n}_{S,S})^{-1}X_S^TX_{S^c}^{\perp}\Theta_{S^c,j}\right|.
\end{align*}
By the Gaussian property of $X$, $X_{S^c}^{\perp}$ is Gaussian and is independent of $X_S$. Conditioning on $X_S$ and $\eps$, we have $x_{i,S^c}^{\perp}\Theta_{S^c,j}\sim N(0,\Theta_{j,S^c}\Theta^{-1}_{S^c,S^c}\Theta_{S^c,j})$. Notice that
\[
     \Theta_{j,S^c}\Theta_{S^c,S^c}^{-1}\Theta_{S^c,j}\leq \Theta_{j,j}\leq K_1.
\]
Since $\hat{\beta}_S$ is a function of $X_S$ and $\eps$, $X^{\perp}_{S^c}\Theta_{S^c,j}$ is independent of $X_S$ and $\hat{\beta}_S$. Hence,
\[
 |r_{j,2}|=O_P\left(\lam\sqrt{\frac{ sgn(\hat{\beta}_S)^T(\Sig^{n}_{S,S})^{-1}sgn(\hat{\beta}_S)K_1}{n}}\right)=O_P\left(\lam\sqrt{\frac{s}{n}}\right),
\]
where we use the fact that
\[
  sgn(\hat{\beta}_S)^T(\Sig^{n}_{S,S})^{-1}sgn(\hat{\beta}_S)\leq \|sgn(\hat{\beta}_S)\|_2^2\|(\Sig^{n}_{S,S})^{-1}\|_2=O_P(s/C_{\min}).
\]
Together with the bound for $|r_{j,1}|$, for $j\notin S$ and $\lam\asymp\sqrt{\log p/n}$, we arrive at
\[
  |\widehat{rem}_j^o|\leq |r_{j,1}|+|r_{j,2}|=O_P\left(\frac{\sqrt{s\log p}}{n}\right).
\]

(ii) If $j\in S$, we use an leave-one-out argument. 
In the event that $\widehat{S}\subseteq S$, for $\hat{\beta}^{(k)}$ defined via (\ref{beta-loo}), it holds that
\begin{align*}
\hat{\beta}^{(DB)}_j-\beta_j&=(e_j^T-\Theta_{j}^T\Sig^n)\hat{u}+\frac{\langle \hat{z}_j,\eps\rangle}{n}\\
&=\underbrace{(e_j^T-\Theta_{j}^T\Sig^n)_{A_j}(\hat{\beta}^{(j)}_{A_j}-\beta_{A_j})}_{r_{j,1}}+ \underbrace{(e_j^T-\Theta_{j}^T\Sig^n)_S(\hat{\beta}_S-\hat{\beta}^{(j)}_S)}_{r_{j,2}}+ \frac{\langle \hat{z}_j,\eps\rangle}{n},
\end{align*}
where the last step is due to $\hat{\beta}^{(j)}_j=\beta_j$.
For $r_{j,1}$, we can bound it similarly as in the case $j\notin S$. In fact, the KKT conditions regarding $\hat{\beta}^{(j)}_{A_j}$ gives
\begin{align}
\label{eq-loo}
   &\Sig^n_{A_j,A_j}(\hat{\beta}^{(j)}_{A_j}-\beta_{A_j})-W^n_{A_j}=-\lam sgn(\hat{\beta}^{(j)}_{A_j}),
\end{align}
where $\Sig^n_{A_j,A_j}$ is invertible with high probability. As a result,
\begin{align*}
r_{j,1}&=-\Theta_{j}^T\Sig^n_{.,A_j}(\Sig^n_{A_j,A_j})^{-1}W^n_{A_j}+\Theta_{j}^T\Sig^n_{.,A_j}(\Sig^n_{A_j,A_j})^{-1}\lam sgn(\hat{\beta}^{(j)}_{A_j})\\
&=-\Theta_{j}^T\Sig^n_{.,A_j}(\Sig^n_{A_j,A_j})^{-1}W^n_{A_j}+\frac{\lam}{n}\Theta_{j,A_j^c}\langle X_{A_j^c}^{\perp},X_{A_j}\rangle(\Sig^n_{A_j,A_j})^{-1}sgn(\hat{\beta}^{(j)}_{A_j}),
\end{align*}
where the last step is due to $\Theta_{A_j,k}+\Sig^{-1}_{A_j,A_j}\Sig_{A_j,A_j^c}\Theta_{A_j^c,k}=0$ for $k\in A_j^c$.
By similar arguments as for the case $j\in S^c$, we have
\[
  |r_{j,1}|=O_P\left(\frac{\sqrt{s\log p}}{n}\right).
\]
For $r_{j,2}$, we have
\begin{align*}
\left|(e_j^T-\Theta_{j}^T\Sig^n)_S(\hat{\beta}_S-\hat{\beta}^{(j)}_S)\right|&\leq |(1-\Theta_{j}^T\Sig^n_{.,j})\hat{u}_j|+ \left|(e_j^T-\Theta_{j}^T\Sig^n)_{A_j}(\hat{\beta}_{A_j}-\hat{\beta}^{(j)}_{A_j})\right|\\
&=O_P(\sqrt{K_1/n})|\hat{u}_j| + \|\Theta_{j}^T\Sig^n_{.,A_j}\|_2\|\hat{\beta}_{A_j}-\hat{\beta}^{(j)}_{A_j}\|_2\\
&= O_P(\sqrt{K_1/n})|\hat{u}_j| +O_P(\sqrt{s/n}) \|\hat{\beta}_{A_j}-\hat{\beta}^{(j)}_{A_j}\|_2,
\end{align*}
where last step is due to the independence between $X\Theta_j$ and $X_{A_j}$.
By (\ref{eq3-2}),
\[
   |\hat{u}_j|\leq \|\hat{u}_S\|_2=O_P\left(\sqrt{\frac{s\log p}{n}}\right).
\]
Using Lemma \ref{lem-loo}, $
 \|\hat{\beta}_{A_j}-\hat{\beta}^{(j)}_{A_j}\|_2=O_P(\rho_S\lam).$
Together we have for $j\in S$,
\[
  |\widehat{rem}_j^o|\leq |r_{j,1}|+|r_{j,2}|=O_P\left(\frac{(\rho_S+1)\sqrt{s\log p}}{n}\right).
\]

It is left to deal with the case where \underline{\textbf{$\Theta$ is unknown}}. By definition,
\begin{align*}
\hat{\beta}^{(DB)}_j-\beta_j&=(e_j^T-\frac{\hat{\gam}_j^T\Sig^n}{\langle X_j,X\hat{\gam}\rangle})\hat{u}+\frac{\hat{z}_j^T\eps}{n}.
\end{align*}
Using KKT condition of (\ref{hgam0}), we have
\begin{align*}
 |\widehat{rem}_j^o|= \big|(e_j^T-\frac{\hat{\gam}_j^T\Sig^n}{\langle X_j,X\hat{\gam}\rangle})\hat{u}\big|&\leq \left\|e_j^T-\frac{\hat{\gam}_j^T\Sig^n}{\langle X_j,X\hat{\gam}\rangle}\right\|_{\infty}\|\hat{u}_S\|_1=O_P(s\lam\lam_j).
\end{align*}
\end{proof}

\subsection{Proof of Lemma \ref{lem2b}}
\label{pf-lem2b}
\begin{proof}
When the results of Lemma \ref{lem1b} hold, $\widehat{A}_j\subseteq A_j$.
Let $\widehat{B}_j=A_j\cap  \widehat{A}_j^c$.

(i) We first study the case where \underline{\textbf{$\Theta$ is known}}.
By (\ref{beta-db2}) and (\ref{eq3-2}), in the event $\Omega_0$ (\ref{omega0}) we have
\begin{align}
\hat{\beta}^{(mDB)}_j-\beta_j&=(e^T_j-\hat{w}_j^TX/n)\hat{u}+\frac{\hat{w}_j^T\eps}{n}=-\hat{w}_j^TX_{A_j}\hat{u}_{A_j}/n+\frac{\hat{w}_j^T\eps}{n}\nonumber\\
&=\hat{b}_j+\underbrace{\frac{1}{n}\hat{w}_j^TX_{A_j}\hat{u}_{A_j}-\hat{b}_j}_{\hat{r}_j}+\frac{\hat{w}_j^T\eps}{n},\label{eq1-pf3}
\end{align}
where
\begin{equation}
\label{eq-bj1}
  \hat{b}_j=\frac{1}{n}\hat{w}_j^TX_{A_j}(\Sig^{n}_{A_j,A_j})^{-1}\lam sgn(\hat{\beta}^{(j)}_{A_j})
\end{equation}
for $\hat{\beta}^{(j)}$ defined in (\ref{beta-loo}).
For $\hat{r}_j$, using the definition of $\hat{w}_j$ in (\ref{eq-z2}), we have
\begin{align*}
|\hat{r}_j|
&\leq\underbrace{\left|\frac{\Theta_j^TX^TP_{\widehat{A}_j}^{\perp}X_{A_j}(\hat{\beta}_{A_j}-\hat{\beta}^{(j)}_{A_j})}{\langle X_j, P_{\widehat{A}_j}^{\perp}X\Theta_j\rangle}\right|}_{\hat{r}_{j,1}}+\underbrace{\left|\frac{\Theta_j^TX^TP_{\widehat{A}_j}^{\perp}X_{A_j}(\hat{\beta}^{(j)}_{A_j}-\beta_{A_j})}{\langle X_j, P_{\widehat{A}_j}^{\perp}X\Theta_j\rangle}-\hat{b}_j\right|}_{\hat{r}_{j,2}}
\end{align*}
For the denominators in $\hat{r}_{j,1}$ and $\hat{r}_{j,2}$, noticing that
\begin{align}
\frac{1}{n}|\langle X_j, P_{\widehat{A}_j}^{\perp}X\Theta_j\rangle|&\geq \frac{1}{n}|\langle X^{\perp}_j, P_{\widehat{A}_j}^{\perp}X\Theta_j\rangle|-\frac{1}{n}|\langle P_{\widehat{A}_j}^{\perp}X_{-j}\Sig^{-1}_{-j,-j}\Sig_{-j,j}, X\Theta_j\rangle|\nonumber\\
&=\frac{1}{n\Theta_{j,j}}\|P_{\widehat{A}_j}^{\perp}X\Theta_j\|_2^2-\frac{1}{n\Theta_{j,j}}|\langle P_{\widehat{A}_j}^{\perp}X_{-j}\Theta_{j,-j}, X\Theta_j\rangle|\nonumber\\
&=1-O_P(|\widehat{A}_j|/n+n^{-1/2}),\label{lem2b-eq2}
\end{align}
where the second step is due to $X_j^{\perp}=X\Theta_j$ and last step is due to $P_{\widehat{A}_j}^{\perp}X_{-j}$ is a function of $X_{-j}$ and is independent of $X\Theta_j$.
For the numerator of $\hat{r}_{j,1}$,
\begin{align*}
|\Theta_j^TX^TP_{\widehat{A}_j}^{\perp}X_{A_j}(\hat{\beta}_{A_j}-\hat{\beta}^{(j)}_{A_j})|&\leq \|\Theta_j^TX^TP_{\widehat{A}_j}^{\perp}X_{\widehat{B}_j}\|_2\|\hat{\beta}_{\widehat{B}_j}-\hat{\beta}^{(j)}_{\widehat{B}_j}\|_2\\
&\leq \|\Theta_j^TX^TP_{\widehat{A}_j}^{\perp}X_{\widehat{B}_j}\|_{2}O_P\left(\rho_S\lam\mathbbm{1}(j\in S)\right),
\end{align*}
where the first step is due to the projection and the second step is due to Lemma \ref{lem-loo} if $j\in S$ and $\hat{\beta}_S=\hat{\beta}^{(j)}_S$ if $j\notin S$.
Since $P_{\widehat{A}_j}^{\perp}X_{\widehat{B}_j}$ is a function of $X_{-j}$, which is independent of $X\Theta_j$, we have
\[
    \|\Theta_j^TX^TP_{\widehat{A}_j}^{\perp}X_{\widehat{B}_j}\|_{2}^2=O_P(\Theta_{j,j}\textrm{Tr}(X_{\widehat{B}_j}^TP_{\widehat{A}_j}^{\perp}X_{\widehat{B}_j}))=O_P(|\widehat{B}_j|n).
\]
Since $|\widehat{B}_j|\leq s_-$ with high probability,
we have
\[
  | \hat{r}_{j,1}|=O_P(\rho_S\sqrt{s_-}\lam\mathbbm{1}(j\in S)/\sqrt{n}).
\]
For $\hat{r}_{j,2}$, by the definition of $\hat{\beta}^{(j)}_{A_j}$ in (\ref{beta-loo}), it holds that
\[
 \hat{\beta}_{A_j}^{(j)}-\beta_{A_j}=(\Sig^n_{A_j,A_j})^{-1}X_{A_j}^T\eps/n-\lam(\Sig^n_{A_j,A_j})^{-1}sgn(\hat{\beta}_{A_j}^{(j)})
\]
and hence
\[
  |\hat{r}_{j,2}|=\left|\frac{\Theta_j^TX^TP_{\widehat{A}_j}^{\perp}X_{A_j}(\Sig^n_{A_j,A_j})^{-1}W^n_{A_j}}{\langle X_j, P_{\widehat{A}_j}^{\perp}X\Theta_j\rangle}\right|.
\]
Using the sub-Gaussian property of $\eps$ conditioning on $X$, we have 
\[
   |\hat{r}_{j,2}|=O_P(\sig\|P_{A_j}P^{\perp}_{\widehat{A}_j}X\Theta_j\|_2/n)=O_P(\sqrt{s_-}/n).
\]
where the last step is due to the independence between $X\Theta_j$ and $P^{\perp}_{\widehat{A}_j}X_{A_j}$ and 
 \[
    \|P_{A_j}P^{\perp}_{\widehat{A}_j}X\Theta_j\|^2_2=O_P(\textrm{Tr}(P_{A_j}P^{\perp}_{\widehat{A}_j}))\leq s_-.
 \]
Together we have,
\[
  |\hat{r}_j|=O_P(\rho_S\sqrt{s_-}\lam\mathbbm{1}(j\in S)/\sqrt{n})+o_P(\sqrt{s_-}/n).
\]

For the bootstrapped estimate of bias, $\hat{b}^*_j$ defined in (\ref{eq-hbj}), for $\hat{b}_j$ defined in (\ref{eq-bj1})
\begin{align}
\hat{b}^*_j&=-\hat{w}_j^TX_{A_j}\hat{u}^*_{A_j}/n=\hat{b}_j-\underbrace{\left(\frac{1}{n} \hat{w}_j^TX_{A_j}\hat{u}^*_{A_j}+\hat{b}_j\right)}_{\hat{r}^*_j}.\label{eq2-pf3}
\end{align}
Notice that
 \begin{align*}
|\hat{r}^*_j|&\leq \underbrace{\frac{1}{n} |\hat{w}_j^TX_{A_j}(\hat{\beta}^*_{A_j}-\hat{\beta}^{*,(j)}_{A_j})|}_{\hat{r}^*_{j,1}}+\underbrace{\frac{1}{n} |\hat{w}_j^TX_{A_j}(\hat{\beta}^{*,(j)}_{A_j}-\hat{\beta}^{(j)}_{A_j})-\hat{b}_j|}_{\hat{r}^*_{j,2}}\\
&\quad +\underbrace{\frac{1}{n} |\hat{w}_j^TX_{A_j}(\hat{\beta}^{(j)}_{A_j}-\hat{\beta}_{A_j})|}_{\hat{r}^*_{j,3}}.
\end{align*}
where $\hat{r}^*_{j,3}=\hat{r}_{j,1}$. By Lemma \ref{lem-loo2}, $\hat{r}^*_{j,1}$ can be similarly bounded as for $\hat{r}_{j,1}$. Therefore,
\[
  |\hat{r}^*_{j,1}|+|\hat{r}^*_{j,3}|=O_P(\sqrt{s_-}\rho_S\lam\mathbbm{1}(j\in S)/\sqrt{n}).
\]
For $\hat{r}^*_{j,2}$, we have
{\small
\begin{align*}
|\hat{r}^*_{j,2}|&=|\frac{1}{n}\hat{w}_j^TX_{A_j}(\hat{\beta}^{*,(j)}_{A_j}-\hat{\beta}^{(j)}_{A_j})+\frac{1}{n}\lam \hat{w}_j^TX_{A_j}(\Sig^n_{A_j,A_j})^{-1}sgn(\hat{\beta}^{(j)}_{A_j})|\\
&\leq |\frac{\lam}{n\lam}\hat{w}_j^TX_{A_j}(\Sig^n_{A_j,A_j})^{-1}W^*_{A_j}|+ \underbrace{\frac{\lam }{n}|\hat{w}_j^TX_{A_j}(\Sig^n_{A_j,A_j})^{-1}[sgn(\hat{\beta}^{(j)}_{A_j})-sgn(\hat{\beta}^{*,(j)}_{A_j})]|}_{\tilde{r}^*_{j,2}},
\end{align*}
}
where the first term can be similarly bounded as $\hat{r}_{j,2}$. The second term can be rewritten using the definition of $\hat{w}_j$:
\begin{align*}
\tilde{r}^*_{j,2}=\frac{\lam }{n}\frac{|\Theta_j^TX^TP_{\widehat{A}_j}^{\perp}X_{A_j}(\Sig^n_{A_j,A_j})^{-1}[sgn(\hat{\beta}^{(j)}_{A_j})-sgn(\hat{\beta}^{*,(j)}_{A_j})]|}{|X_j^TP_{\widehat{A}_j}^{\perp}X\Theta_j|}
\end{align*}
Noticing that  $\hat{\beta}^{*,(j)}_{A_j}$, and $\hat{\beta}_{A_j}^{(j)}$ are functions of $X_{-j}$ and $\eps$. Hence, $X\Theta_j$ is independent of $P_{\widehat{A}_j}^{\perp}X_{A_j}$, $\hat{\beta}^{*,(j)}_{A_j}$, and $\hat{\beta}_{A_j}^{(j)}$ no matter $j\in S$ or not. Hence,
\begin{align*}
  \tilde{r}^*_{j,2}&=O_P\left(\frac{\lam}{n}\|P_{\widehat{A}_j}^{\perp}X_{A_j}(\Sig^n_{A_j,A_j})^{-1}[sgn(\hat{\beta}^{(j)}_{A_j})-sgn(\hat{\beta}^{*,(j)}_{A_j})]\|_2\right)\\
   &=O_P\left(\frac{\lam}{\sqrt{n}}\|(\Sig^n_{A_j,A_j})^{-1}\|^{1/2}_2\|sgn(\hat{\beta}^{(j)}_{A_j})-sgn(\hat{\beta}^{*,(j)}_{A_j})\|_2\right).
\end{align*}
If $j\notin S$, $\hat{\beta}^{(j)}_S=\check{b}_S$ defined in (\ref{cbs}) and $\hat{\beta}^{*,(j)}_S=\check{b}^*_S$ defined in (\ref{cbs-star}). We have
\begin{align*}
  &\|sgn(\hat{\beta}^{(j)}_{A_j})-sgn(\hat{\beta}^{*,(j)}_{A_j})\|_2\\
  &\leq\|sgn(\hat{b}^*_{A_j})-sgn(\beta_{A_j})\|_2+\|sgn(\hat{\beta}^{*}_{A_j})-sgn(\beta_{A_j})\|_2\leq \sqrt{s_-}.
\end{align*}
If $j\in S$, we have for a large enough constant $C$, it holds with high probability that
\[
   \|\hat{\beta}^{(j)}_{A_j}-\hat{\beta}^{*,(j)}_{A_j}\|_{\infty}\leq \|\hat{\beta}^{(j)}_{A_j}-\hat{\beta}_{A_j}\|_{\infty}+\|\hat{\beta}^*_{A_j}-\hat{\beta}^{*,(j)}_{A_j}\|_{\infty}+\|\hat{\beta}^*_{A_j}-\hat{\beta}_{A_j}\|_{\infty}\leq C\rho_S\lam,
\]
where the last step is due to Lemma \ref{lem-loo}, Lemma \ref{lem-loo2} and Lemma \ref{lem1b}.
Therefore, with high probability,
 \begin{align*}
    & \|sgn(\hat{\beta}^{(j)}_{A_j})-sgn(\hat{\beta}^{*,(j)}_{A_j})\|^2_2\leq |\{k\in A_j: |\hat{\beta}^{(j)}_{k}|\leq C^*\rho_S\lam \}|\\
    &\leq  |\{k\in A_j: |\beta_{k}|\leq (C+C^*)\rho_S\lam \}|\leq |S_-|.
 \end{align*}
Hence, we have proved
\[
  |\hat{r}^*_j|=O_P(\sqrt{s_-}\lam/\sqrt{n}+\sqrt{s_-}\rho_S\lam/\sqrt{n}\mathbbm{1}(j\in S)).
\]
Invoking (\ref{eq1-pf3}) and (\ref{eq2-pf3}),
\[
   \hat{\beta}^{(BS-DB)}_j-\beta_j=\frac{\hat{w}_j^T\eps}{n}+\hat{b}_j-\hat{r}_j-(\hat{b}_j-\hat{r}^*_j)=\frac{\hat{w}_j^T\eps}{n}+|\hat{r}_j|+|\hat{r}^*_j|.
\]
The proof for known $\Theta$ is complete.

(ii) Next we consider \underline{\textbf{unknown $\Theta$}}.
By (\ref{beta-db2}) and Lemma \ref{lem1a} we have
\begin{align*}
\hat{\beta}^{(mDB)}_j-\beta_j&=(e^T_j-\hat{w}_j^TX/n)\hat{u}+\frac{\hat{w}_j^T\eps}{n}=-\hat{w}_j^TX_{A_j}\hat{u}_{A_j}/n+\frac{\hat{w}_j^T\eps}{n}\\
&=\frac{1}{n}\hat{w}_j^TX_{A_j}e_{A_j}^T(\Sig^{n}_{S,S})^{-1}[\lam sgn(\hat{\beta}_S)-W^{n}_S]+\frac{\hat{w}_j^T\eps}{n}\\
&=\hat{b}_j-\underbrace{\frac{1}{n}\hat{w}_j^TX_{A_j}e_{A_j}^T(\Sig^{n}_{S,S})^{-1}W^{n}_S}_{\hat{r}_j}+\frac{\hat{w}_j^T\eps}{n},
\end{align*}
where 
\begin{equation}
\label{eq-bj2}
 \hat{b}_j=-\frac{1}{n} \hat{w}_j^TX_{A_j}e_{A_j}^T(\Sig^{n}_{S,S})^{-1}\lam sgn(\hat{\beta}_S).
 \end{equation}
 By the KKT condition of (\ref{hkappa}),
\begin{equation}
\label{eq4-2}
  \Sig^n_{\widehat{A}_j,.}\hat{\kappa}=0~\text{and}~ \Sig^n_{\widehat{S}^c\setminus\{j\},.}\hat{\kappa}=-\lam_j sgn(\hat{\kappa}_{\widehat{S}^c\setminus\{j\}}).
\end{equation}
As a result,
  \begin{align*}
  |\hat{r}_j|&=\left|\frac{1}{n}\hat{w}_j^TX_{\widehat{B}_j}e_{\widehat{B}_j}^T(\Sig^{n}_{S,S})^{-1}W^{n}_S\right|\leq |\widehat{B}_j|\|\hat{w}_j^TX_{\widehat{B}_j}/n\|_{\infty}\|e_{\widehat{B}_j}^T(\Sig^{n}_{S,S})^{-1}W^{n}_S\|_{\infty}\\
  &\leq  |\widehat{B}_j|\|\hat{w}_j^TX_{\widehat{B}_j}/n\|_{\infty}\max_{k\in S_-}\|e_{k}^T(\Sig^{n}_{S,S})^{-1}W^{n}_S\|_{\infty}=O_P(s_-\lam_j\lam),
  \end{align*}
  where the first and last step are due to the (\ref{eq4-2}), the third step is due to $\widehat{B}_j\subseteq  S_-$ with high probability. 
  
For the bootstrapped estimate of bias,
\begin{align*}
\hat{b}^*_j&=-\hat{w}_j^TX_{A_j}\hat{u}^*_{A_j}/n=\frac{\lam}{n} \hat{w}_j^TX_{A_j}e_{A_j}^T(\Sig^{n}_{S,S})^{-1} sgn(\hat{\beta}^*_S)\\
&=\hat{b}_j+\underbrace{\frac{1}{n} \hat{w}_j^TX_{A_j}e_{A_j}^T(\Sig^{n}_{S,S})^{-1}[\lam sgn(\hat{\beta}^*_S)-\lam sgn(\hat{\beta}_S)]}_{\hat{r}^*_j}.
\end{align*}
Now we bound $\hat{r}^*_j$. Using (\ref{eq4-2}) again, we have

\begin{align*}
|\hat{r}^*_j|=&\lam|\frac{1}{n} \hat{w}_j^TX_{A_j}e_{A_j}^T(\Sig^{n}_{S,S})^{-1}[sgn(\hat{\beta}^*_S)-sgn(\hat{\beta}_S)]|\\
&=\frac{n\lam\lam_j}{\langle X_j,P_{\widehat{A}_j}^{\perp}X\hat{\kappa}\rangle}\left|sgn(\hat{\kappa}_{\widehat{B}_j})^Te_{\widehat{B}_j}^T(\Sigma^{n}_{S,S})^{-1}e_{S_-}(sgn(\hat{\beta}^*_{S_-}) -sgn(\hat{\beta}_{S_-})\right|\\
&\leq \lam\lam_j(1+o_P(1))\|sgn(\hat{\kappa}_{\widehat{B}_j})\|_2\|(\Sigma^{n}_{S,S})^{-1}\|_2\|sgn(\hat{\beta}^*_{S_-}) -sgn(\hat{\beta}_{S_-})\|_2\\
&\leq O_P(s_-\lam\lam_j/C_{\min}),
\end{align*}

where the first step is by (\ref{eq4-2}), the second step is due to (\ref{lem2b-eq2}) and the last step is due to Lemma \ref{lem1b}.
Therefore, $  \hat{r}^*_j=O_P(s_-\lam\lam_j)$.

\end{proof}

\begin{proof}[Proof of Theorem \ref{thm2b}]
By Lemma \ref{lem2b}, we are left to establish the asymptotic normality of $\langle \hat{w}_j,\eps\rangle /\sqrt{n}$. 

First consider the case where \underline{$\Theta_j$ is known}. Since $P^{\perp}_{\widehat{A}_j}\eps$ is a function of $\eps$ and $X_{-j}$, it is independent of $X\Theta_j$. As $\langle X\Theta_j,P^{\perp}_{\widehat{A}_j}\eps\rangle|X_{-j},\eps \sim N(0,\|P^{\perp}_{\widehat{A}_j}\eps\|_2^2\Theta_{j,j})$ and $\langle X_j,P^{\perp}_{\widehat{A}_j}X\Theta_j\rangle/n =\Theta_{j,j}+o_P(1)$, we have for any $t\in (-\infty,\infty)$,
\begin{align*}
&\P\left(\langle \hat{w}_j,\eps\rangle /\sqrt{n}\leq t|X_{-j},\eps\right)=\Phi(\Theta^{1/2}_{j,j}t/\sig)+o(1).
\end{align*}
Together with Lemma \ref{lem2b},
\[
 \P\left(\sqrt{n}(\hat{\beta}_j^{(BS-DB)}-\beta_j) \leq t\right)=\Phi(\Theta^{1/2}_{j,j}t/\sig)+o(1).
\]
The fact that $\|\hat{w}_j\|_2/\sqrt{n}=\Theta_{j,j}^{1/2}+o_P(1)$ can be similarly shown as above.

Next, we consider the case where \underline{$\Theta$ is unknown}. Let $\kappa^*=\Theta_{j,j}^{-1}\Theta_{j}$ and $\kappa$ be such that
\[
    X_{\widehat{A}_j}^TX\kappa=0~\text{and}~\kappa_{\widehat{A}_j^c}=\kappa^*_{\widehat{A}_j^c}.
\]
Notice that $\kappa_{\widehat{A}_j}=-(\Sig^n_{\widehat{A}_j,\widehat{A}_j})^{-1}\Sig^n_{\widehat{A}_j,\widehat{A}_j^c}\kappa_{\widehat{A}_j^c}=-(\Sig^n_{\widehat{A}_j,\widehat{A}_j})^{-1}\Sig^n_{\widehat{A}_j,\widehat{A}_j^c}\kappa^*_{\widehat{A}_j^c}$ and $\kappa_j=\kappa^*_j=1$.
Let $\hat{v}=\hat{\kappa}-\kappa$. 
For some positive constant $\phi_0$, define
\begin{align}
\label{eq-e3}
 \mathcal{E}_1&=\left\{\|\Sig^n_{\widehat{S}^c\cap\{-j\},.}\kappa\|_{\infty}\leq \lam_j/2,\right.\nonumber\\
 &\left.\quad\quad \inf_{|J|\leq \|\kappa_{\widehat{S}^c\cap\{-j\}}\|_0}\inf_{\|u_{J^c}\|_1\leq 3\|u_J\|_1\neq 0}\frac{\|P_{\widehat{A}_j}^{\perp}X_{\widehat{S}^c\cap\{-j\}}u\|_2^2}{n\|u_J\|_2^2}\geq \phi_0>0\right\}.
\end{align}
We first prove the desired results in the event  $\mathcal{E}_1$. At the end of the proof, we verify that $\mathcal{E}_1$ holds with high probability.

In view of (\ref{hkappa}), we have
\begin{align*}
   \frac{1}{2n}\|X\hat{v}\|_2^2&\leq |\hat{v}^T\Sig^n\kappa|+\lam_j \|\kappa_{\widehat{S}^c\cap\{-j\}}\|_1-\lam_j \|\hat{\kappa}_{\widehat{S}^c\cap\{-j\}}\|_1\\
  &\quad= |\hat{v}_{\widehat{A}_j^c}^T\Sig^n_{\widehat{A}_j^c,.}\kappa|+\lam_j \|\kappa_{\widehat{S}^c\cap\{-j\}}\|_1-\lam_j \|\hat{\kappa}_{\widehat{S}^c\cap\{-j\}}\|_1\\
   &\quad \leq \|\hat{v}_{\widehat{S}^c\cap\{-j\}}\|_1\|\Sig^n_{\widehat{S}^c\cap\{-j\},.}\kappa\|_{\infty}+\lam_j \|\kappa_{\widehat{S}^c\cap\{-j\}}\|_1-\lam_j \|\hat{\kappa}_{\widehat{S}^c\cap\{-j\}}\|_1,
\end{align*}
where the second step follows from the definition of $\kappa$ and the last step is due to $\kappa_j=\hat{\kappa}_j=1$. 
In event $\mathcal{E}_1$, the following oracle inequality holds:
\begin{align}
\label{eq-ora1}
\frac{1}{2}\hat{v}^T\Sig^n\hat{v}\leq \lam_j\|\kappa_{\widehat{S}^c\cap \{-j\}}\|_1- \lam_j\|\hat{\kappa}_{\widehat{S}^c\cap \{-j\}}\|_1+ \frac{\lam_j}{2}\|\hat{v}_{\widehat{S}^c\cap \{-j\}}\|_1,
\end{align}
where
 \[
  \hat{v}^T\Sig^n\hat{v}=\frac{1}{n}\|P_{\widehat{A}_j}^{\perp}X_{\widehat{S}^c\cap\{-j\}}\hat{v}_{\widehat{S}^c\cap\{-j\}}\|_2^2.
\]
Let $J$ denote the support of $\kappa_{\widehat{S}^c\cap\{-j\}}$. Standard decomposition of the right hand side of (\ref{eq-ora1}) leads to
\begin{align*}
\frac{1}{2n}\|P_{\widehat{A}_j}^{\perp}X_{\widehat{S}^c\cap\{-j\}}\hat{v}_{\widehat{S}^c\cap\{-j\}}\|_2^2\leq \frac{3\lam_j}{2}\sum_{k\in J}|\hat{v}_k|- \frac{\lam_j}{2}\sum_{k\in \widehat{S}^c\cap \{-j\}\setminus J}|\hat{v}_k|.
\end{align*}
Using the second statement in event $\mathcal{E}_1$, we have
\[
   \frac{\phi_0}{2}\|\hat{v}_J\|_2^2\leq \frac{3\lam_j}{2}\sum_{k\in J}|\hat{v}_k|\leq \frac{3\lam_j}{2}\sqrt{|J|}\|\hat{v}_J\|_2.
\]
Hence,
\begin{align*}
&|\hat{v}_J\|_2\leq \frac{3\lam_j\sqrt{|J|}}{\phi_0},\quad\|\hat{v}_J\|_1\leq \sqrt{|J|}\|\hat{v}_J\|_2\leq \frac{3\lam_j|J|}{\phi_0},~\text{and}~\\
&\|\hat{v}_{\widehat{S}^c\cap\{-j\}}\|_1\leq (1+3)\|\hat{v}_J\|_1\leq \frac{12\lam_j|J|}{\phi_0}.
\end{align*}
According to our choice of $\lam_j$, we arrive at
\[
   \|\hat{v}_{\widehat{A}_j^c}\|_1=O_P\left(\frac{\|\kappa_{\widehat{S}^c\cap\{-j\}}\|_0}{\phi_0}\sqrt{\frac{\log p}{n}}\right).
\]
Therefore,
\begin{align*}
   &\frac{\langle \eps,X\hat{\kappa}\rangle}{n}= \frac{\langle \eps,X\kappa\rangle}{n}+\frac{\langle \eps,X_{\widehat{S}^c\cap\{-j\}}\hat{v}_{\widehat{S}^c\cap\{-j\}}\rangle}{n}\\
   &=\frac{\langle \eps,X\kappa\rangle}{n}+\left\|\frac{X^T\eps}{n}\right\|_{\infty}\|\|\hat{v}_{\widehat{S}^c\cap\{-j\}}\|_1=\frac{\langle \eps,X\kappa\rangle}{n}+O_P\left(\frac{\|\Theta_{\widehat{S}^c\cap\{-j\},j}\|_0\log p}{n}\right),
\end{align*}
where $\langle \eps,X\kappa\rangle/n=\langle \eps,P_{\widehat{A}_j}^{\perp}X_{\widehat{A}_j^c}\kappa^*_{\widehat{A}_j^c}\rangle/n$. The rest of proof follows from the above statement for known $\Theta$.

It is left to verify $\P(\mathcal{E}_1)\rightarrow 1$ for $\mathcal{E}_1$ defined in (\ref{eq-e3}).  
We first notice that
\begin{align*}
\Sig^n_{\widehat{S}^c\cap\{-j\},.}\kappa&=\frac{1}{n}\langle X_{\widehat{S}^c\cap\{-j\}}, X_j+X_{\widehat{A}_j}\kappa_{\widehat{A}_j}+X_{\widehat{S}^c\cap\{-j\}}\kappa_{\widehat{S}^c\cap\{-j\}}\rangle\\
&=\frac{1}{n}\langle X_{\widehat{S}^c\cap\{-j\}}, P_{\widehat{A}_j}^{\perp}(X_j+X_{\widehat{S}^c\cap\{-j\}}\kappa_{\widehat{S}^c\cap\{-j\}})\rangle\\
&=\frac{1}{n}\langle X_{\widehat{S}^c\cap\{-j\}}, P_{\widehat{A}_j}^{\perp}X_{\widehat{A}_j^c}\kappa_{\widehat{A}_j^c}\rangle=\frac{1}{n}\langle X_{\widehat{S}^c\cap\{-j\}}, P_{\widehat{A}_j}^{\perp}X\kappa^*\rangle,
\end{align*}
where the last step is due to $P_{\widehat{A}_j}^{\perp}X_{\widehat{A}_j^c}\kappa_{\widehat{A}_j^c}=P_{\widehat{A}_j}^{\perp}X_{\widehat{A}_j^c}\kappa_{\widehat{A}_j^c}+P_{\widehat{A}_j}^{\perp}X_{\widehat{A}_j}\kappa^*_{\widehat{A}_j}$ and $\kappa_{\widehat{A}_j^c}=\kappa^*_{\widehat{A}_j^c}$.
By the Gaussian property of $X$, $X\kappa^*$ is independent of $X_{-j}$. As $P_{\widehat{A}_j}^{\perp}X_{\widehat{S}^c\cap\{-j\}}$ is a function of $X_{-j}$,
\[
  \P\left(\|\Sig^n_{\widehat{S}^c\cap\{-j\},.}\kappa\|_{\infty}\geq c_1\sqrt{\log p/n}\right)=o(1)
\]
for some positive constant $c_1$.
Hence, for $\lam_j=C\sqrt{\log p/n_0}$ with large enough constant $C$, the first statement in event $\mathcal{E}_1$ holds with high probability.

For the second statement, in the event that $S_+\subseteq \tilde{S}\subseteq S$, it holds that
{\small
\begin{align*}
&\P\left(\inf_{|J|\leq \|\kappa_{\widehat{S}^c\cap\{-j\}}\|_0}\inf_{\|u_{J^c}\|_1\leq 3\|u_J\|_1\neq 0}\frac{\|P_{\widehat{A}_j}^{\perp}X_{\widehat{S}^c\cap\{-j\}}u\|_2^2}{n\|u_J\|_2^2}\leq \phi_0\right)\\
&\leq 2^{s_{-}}\max_{S_+\subseteq\tilde{S}\subseteq S}\P\left(\inf_{|J|\leq \|\kappa_{\tilde{S}^c\cap\{-j\}}\|_0}\inf_{\|u_{J^c}\|_1\leq 3\|u_J\|_1\neq 0}\frac{\|P_{\tilde{A}_j}^{\perp}X_{\tilde{S}^c\cap\{-j\}}u\|_2^2}{n\|u_J\|_2^2}\leq \phi_0\right)+o(1),
\end{align*}
}
where $\tilde{A}_j=\tilde{S}\cap\{-j\}$ and inequality is because there are at most $2^{s_-}$ different $\tilde{S}$.
We further use the fact that
\begin{align*}
&\inf_{\|u_{J^c}\|_1\leq 3\|u_J\|_1\neq 0}\frac{\|P_{\tilde{A}_j}^{\perp}X_{\tilde{S}^c\cap\{-j\}}u\|_2^2}{n\|u_J\|_2^2}= \inf_{\|u_{J^c}\|_1\leq 3\|u_J\|_1\neq 0}\frac{\||P_{\tilde{A}_j}^{\perp}X^{\perp}_{\tilde{S}^c\cap\{-j\}}u\|_2^2}{n\|u_J\|_2^2},
\end{align*}
where $X^{\perp}_{\tilde{S}^c\cap\{-j\}}=X_{\tilde{S}^c\cap\{-j\}}-X_{\tilde{A}_j}\Sig_{\tilde{A}_j,\tilde{A}_j}^{-1}\Sig_{\tilde{A}_j,\tilde{S}^c\cap\{-j\}}$
is independent of $X_{\tilde{A}_j}$. Let $\{g_1,\dots,g_{n-|\tilde{A}_j|}\}$ be an orthonormal basis for the complement of the column space of $X_{\tilde{A}_j}$.
Then we have 
\begin{align*}
\|P_{\tilde{A}_j}^{\perp}X^{\perp}_{\tilde{S}^c\cap\{-j\}}u\|_2^2&=\sum_{k=1}^{n-|\tilde{A}_j|}\|g_k^TX^{\perp}_{\tilde{S}^c\cap\{-j\}}u\|_2^2,
\end{align*}
where $g_k^TX^{\perp}_{\tilde{S}^c\cap\{-j\}}|X_{\tilde{A}_j}$ is Gaussian with mean zero and variance $g_k^T\Sig^{\perp}_{\tilde{S}^c\cap\{-j\},\tilde{S}^c\cap\{-j\}}g_k$, where
\[
   \Sig^{\perp}_{\tilde{S}^c\cap\{-j\},\tilde{S}^c\cap\{-j\}}=\Sig_{\tilde{S}^c\cap\{-j\},\tilde{S}^c\cap\{-j\}}-\Sig_{\tilde{S}^c\cap\{-j\},\tilde{A}_j}\Sig_{\tilde{A}_j,\tilde{A}_j}^{-1}\Sig_{\tilde{A}_j,\tilde{S}^c\cap\{-j\}}.
\]
Moreover, $v_k^TX^{\perp}_{\tilde{S}^c\cap\{-j\}}$ is independent of $v_{k'}^TX^{\perp}_{\tilde{S}^c\cap\{-j\}}$ for $k\neq k'$ conditioning on $X_{\tilde{A}_j}$.
 We can then use Theorem 1.6 of \citet{Zhou09} to show that
 \begin{align*}
&2^{s_-} \P\left(\inf_{|J|\leq \|\kappa_{\tilde{S}^c\cap\{-j\}}\|_0}\inf_{\|u_{J^c}\|_1\leq 3\|u_J\|_1\neq 0}\sum_{k=1}^{n-|\tilde{A}_j|}\|g_k^TX^{\perp}_{\tilde{S}^c\cap\{-j\}}u\|_2^2\leq \phi_0\right)\\
 &\leq 2^{s_-}\P\left(\inf_{|J|\leq \|\kappa_{(S_+)^c}\|_0}\inf_{\|u_{J^c}\|_1\leq 3\|u_J\|_1\neq 0}\sum_{k=1}^{n-s}\|g_k^TX^{\perp}_{\tilde{S}^c\cap\{-j\}}u\|_2^2\leq \phi_0\right)=o(1)
 \end{align*}
 if $n-s\gg \max\{\max\{\|\kappa_{(S_+)^c}\|_0,1\}\log p , s_-\}$ and $\Lambda_{\min}(\Sig^{\perp}_{\tilde{S}^c\cap\{-j\},\tilde{S}^c\cap\{-j\}})\geq 2\phi_0$. In fact,
 \[
    \{\Sig^{\perp}_{\tilde{S}^c\cap\{-j\},\tilde{S}^c\cap\{-j\}}\}^{-1}=\{(\Sig_{-j,-j})^{-1}\}_{\tilde{S}^c\cap\{-j\},\tilde{S}^c\cap\{-j\}},
 \]
 which is positive definite for any $S_+\subseteq \tilde{S} \subseteq S$.
\end{proof}

\section*{Acknowledgements}
The author would like to thank Professor Cun-Hui Zhang for insightful discussions on the technical part of the paper as well as the presentation.
\bibliography{BootBib}{}
\bibliographystyle{chicago}

\appendix


\section{Proof of lemmas and theorems in Section 3}
\subsection{Some technical lemmas}

\begin{lemma}
\label{alem2}
Under Condition \ref{cond1b} - Condition \ref{cond3b}, we have the following results.
Let $c_1>4$, $c_2>0$. For $n>c_1s$, with probability at least $1-2\exp(-c_2n)$,
\begin{equation}
\label{eq: lem4a}
  C_{\min}/4\leq \Lambda_{\min}(\Sig^n_{S,S})\leq \Lambda_{\max}(({\Sig}^n_{S,S}))\leq 4C_{\max}.
\end{equation}
\end{lemma}
\begin{proof}[Proof of Lemma \ref{alem2}]
The proof is standard. See Corollary 5.35 of \citep{Ver10}.
\end{proof}

\begin{lemma}Under conditions \ref{cond1b} - \ref{cond3b}, if $s\log p\leq c_0n$ for a small enough constant $c_0$, then it holds that with with probability at least $1-c_1/p-c_2s/n-\exp(-c_3n)$, there exists a large enough constant $C$ such that for $\hat{\beta}^{(k)}_{A_k}$ defined in (\ref{beta-loo}),
\label{lem-loo}
 \[
     \max_{k\in S}\left\|\hat{\beta}_{A_k}-\hat{\beta}^{(k)}_{A_k}\right\|_2\leq C\rho_S\lam.
 \]
\end{lemma}
\begin{proof}[Proof of Lemma \ref{lem-loo}]
We justify Lemma \ref{lem-loo} in the event $\{C_1\leq \Lambda_{\min}(\Sig^n_{S,S})\leq \Lambda_{\max}(\Sig^n_{S,S})\leq C_2\}$, which holds with high probability by Lemma \ref{alem2}.
We similarly show the results of Lemma 6.4 - Lemma 6.7 in \citet{JM15} under current conditions. 
Specifically, define
\[
    u_k(b)=\frac{X_k^T(\eps+X_{A_k}(\beta_{A_k} - b))}{n}.
\]
%
By (82) of \citet{JM15}, we have
\begin{align}
\label{eq-ap2}
\frac{1}{2n}\left\|P_{X_k}^{\perp}X_{A_k}(\hat{\beta}_{A_k}-\hat{\beta}^{(k)}_{A_k})\right\|_2^2\leq |u(\hat{\beta}^{(k)}_{A_k})|\left|u(\hat{\beta}_{A_k})-u(\hat{\beta}^{(k)}_{A_k})\right|.
\end{align}
Noticing that
 \begin{align}
 \max_{k\in S}|u(\hat{\beta}^{(k)}_{A_k})|&= \max_{k\in S}\left|\frac{X_k^T(\eps+X_{A_k}(\beta_{A_k} - \hat{\beta}^{(k)}_{A_k}))}{n}\right|\nonumber\\
 &\leq  \max_{k\in S}\left|\frac{X_k^T\eps}{n}\right|+ \max_{k\in S}\left|\Sig_{k,A_k}\Sig_{A_k,A_k}^{-1}\Sig^n_{A_k,A_k}(\beta_{A_k} - \hat{\beta}^{(k)}_{A_k})\right|\nonumber\\
 &+ \max_{k\in S}\left|\frac{\langle X_k^{\perp},X_{A_k}(\beta_{A_k} - \hat{\beta}^{(k)}_{A_k})\rangle}{n}\right|\nonumber,
 \end{align}
 where the first term is no larger that $c_1\sqrt{\log p/n}$ with probability at least $1-c_2/p$, the second term can be bounded by
 \begin{align*}
  & \max_{k\in S}\frac{1}{e_k^T\Sig_{S,S}^{-1}e_k} |e_k^T\Sig_{S,S}^{-1}e_{A_k}\Sig^n_{A_k,A_k}(\beta_{A_k} - \hat{\beta}^{(k)}_{A_k})|\\
   &\leq  \frac{1}{C_{\min}}\rho_S\max_{k\in S} \|\Sig^n_{A_k,A_k}(\beta_{A_k} - \hat{\beta}^{(k)}_{A_k})\|_{\infty}  \leq  \frac{\lam \rho_S}{C_{\min}},
    \end{align*}
    where the last step is due to the KKT condition of $\hat{\beta}^{(k)}_{A_k}$.
    With probability at least $1-c_3/p$, it holds that
    \[
    \max_{k\in S}\left|\frac{\langle X_k^{\perp},X_{A_k}(\beta_{A_k} - \hat{\beta}^{(k)}_{A_k})\rangle}{n}\right|\leq c_4K_1\sqrt{\log p}/n\max_{k\in S}\|X_{A_k}(\beta_{A_k} - \hat{\beta}^{(k)}_{A_k})\|_2.
    \]
It is left to show that
\[
   \frac{1}{\sqrt{n}} \max_{k\in S}\|X_{A_k}(\beta_{A_k} - \hat{\beta}^{(k)}_{A_k})\|_2=o_P(1).
\]
Using the KKT condition of (\ref{beta-loo}),
\[
   \hat{\beta}^{(k)}_{A_k}-\beta_{A_k}=(\Sig^n_{A_k,A_k})^{-1}W^n_{A_k}-\lam (\Sig^n_{A_k,A_k})^{-1}sgn(\hat{\beta}^{(k)}_{A_k})
\]
and hence
\[
   \frac{1}{n}\|X_{A_k}( \hat{\beta}^{(k)}_{A_k}-\beta_{A_k})\|_2^2\leq \frac{2}{n}\|P_{A_k}\eps\|_2^2+2\lam^2sgn(\hat{\beta}^{(k)}_{A_k})^T(\Sig^n_{A_k,A_k})^{-1}sgn(\hat{\beta}^{(k)}_{A_k}),
\]
where $P_{A_k}=X_{A_k}(X_{A_k}^TX_{A_k})^{-1}X^T_{A_k}$. Using the fact that $A_k\subseteq S$ for any $k\in S$, we further have
\begin{align}
&\max_{k\in S} \frac{1}{n}\|X_{A_k}( \hat{\beta}^{(k)}_{A_k}-\beta_{A_k})\|_2^2\nonumber\\
&\leq \max_{k\in S}\frac{2}{n}\|P_{A_k}\eps\|_2^2 + 2\lam^2\max_{k\in S}sgn(\hat{\beta}^{(k)}_{A_k})^T(\Sig^n_{A_k,A_k})^{-1}sgn(\hat{\beta}^{(k)}_{A_k})\nonumber\\
&\leq \frac{2}{n}\|P_S\eps\|_2^2+ 2s\lam^2 \max_{k\in S}\|(\Sig^n_{A_k,A_k})^{-1}\|_2\leq \frac{c_4K_2^2s}{n}+ \frac{4s\lam^2}{C_{\min}},\label{eq-tlem2-5}
\end{align}
with probability at least $1-c_5s/n-\exp(-c_6n)$.
The last step is by Cheybeshev's inequality and the fact that
\[
   (\Sig^n_{A_k,A_k})^{-1}\preceq (\Sig^n_{S,S})_{A_k,A_k}^{-1}~\text{and}~\|(\Sig^n_{S,S})_{A_k,A_k}^{-1}\|_2\leq \|(\Sig^n_{S,S})^{-1}\|_2.
\]

We arrive at if $s\log p/n\leq c_7$ for small enough constant $c_7$, then with probability at least $1-c_3/p-c_5s/n-\exp(-c_6n)$, there  exists a large enough constant C such that
\begin{equation}
 \label{eq-ap3}
 \max_{k\in S}|u(\hat{\beta}^{(k)}_{A_k})|\leq C\rho_S\lam.
\end{equation}

 Note that
 \[
 \max_{k\in S}\frac{\left|u(\hat{\beta}_{A_k})-u(\hat{\beta}^{(k)}_{A_k})\right|}{\left\|\hat{\beta}_{A_k}-\hat{\beta}^{(k)}_{A_k}\right\|_2}\leq \max_{k\in S}\left\|\frac{X_k^TX_{A_k}}{n}\right\|_2\leq C_2.
 \]
 Together with (\ref{eq-ap2}) and (\ref{eq-ap3}), we have
 \begin{align}
\max_{k\in S} \frac{\left\|P_{k}^{\perp}X_{A_k}(\hat{\beta}_{A_k}-\hat{\beta}^{(k)}_{A_k})\right\|_2^2}{2n\left\|\hat{\beta}_{A_k}-\hat{\beta}^{(k)}_{A_k}\right\|_2}&\leq  \max_{k\in S}|u(\hat{\beta}^{(k)}_{A_k})| \max_{k\in S}\frac{\left|u(\hat{\beta}_{A_k})-u(\hat{\beta}^{(k)}_{A_k})\right|}{\left\|\hat{\beta}_{A_k}-\hat{\beta}^{(k)}_{A_k}\right\|_2}\nonumber\\
 &\leq C_2C\rho_S\lam\label{eq-ap4}
 \end{align}
 with high probability.
 Define $\mu^{(k)}\in\R^s$ such that
 \[
   \mu^{(k)}_k=-\left\langle \frac{X_k}{\|X_k\|^2_2},X_{A_k}(\hat{\beta}_{A_k}-\hat{\beta}^{(k)}_{A_k})\right\rangle,~~\mu^{(k)}_{A_k}= \hat{\beta}_{A_k} - \hat{\beta}^{(k)}_{A_k}.
 \]
 \begin{align*}
 \max_{k\in S} \frac{\left\|P_{k}^{\perp}X_{A_k}(\hat{\beta}_{A_k}-\hat{\beta}^{(k)}_{A_k})\right\|_2^2}{2n\left\|\hat{\beta}_{A_k}-\hat{\beta}^{(k)}_{A_k}\right\|^2_2}&=\max_{k\in S} \frac{\|X_S\mu^{(k)}\|_2^2}{2n\left\|\hat{\beta}_{A_k}-\hat{\beta}^{(k)}_{A_k}\right\|^2_2}
 \geq \frac{C_1}{2}.
 \end{align*}
 Together with (\ref{eq-ap4}), we have
 \[
     \left\|\hat{\beta}_{A_k}-\hat{\beta}^{(k)}_{A_k}\right\|_2\leq \frac{2C_2C\rho_S\lam}{C_1}.
 \]
with probability at least $1-c_3/p-c_5s/n-\exp(-c_6n)$. 
\end{proof}

\subsection{Proof of Lemma \ref{lem1a}}
\label{pf-lem1a}

\begin{proof}[Proof of Lemma \ref{lem1a}]
(i)
 First consider a restricted Lasso problem
\begin{equation}
\label{cbs}
   \check{b}_S=\argmin_{b\in \R^s} \left\{\frac{1}{2n}\|y-X_Sb_S\|_2^2+\lam\|b_S\|_1\right\}.
\end{equation}
Define
\begin{equation}
\label{eq3-1}
  T^o_k = X_{k}^T\left(X_{S}(X_{S}^TX_{S})^{-1}sgn(\check{b}_S) + \frac{P^{\perp}_S\eps}{n\lam}\right).
\end{equation}
By \citet{Wain09}, if $|T^o_k|<1$ for $\forall k\in S^c$, then Lasso has a unique solution such that $\widehat{S}\subseteq S$.

Since $\check{b}_S$ is only a function of $X_S$ and $\eps$, conditioning on $X_S$ and $\eps$, $T^o_k$ is a Gaussian random variable with mean $\Sig_{k,S}\Sig^{-1}_{S,S}sgn(\check{b}_S)$ and variance  \begin{align*}
 \Var(T^o_k|X_S,\eps)&\leq\Sig_{k,k}\left\|X_{S}(X_{S}^TX_{S})^{-1}sgn(\check{b}_S)+ \frac{P^{\perp}_S\eps}{n\lam}\right\|^2_2\\
 &\leq \Sig_{k,k}( sgn(\check{b}_S)^T(X_{S}^TX_{S})^{-1}sgn(\check{b}_S)+ \|\eps\|_2^2/(n\lam)^2).
 \end{align*}
Let 
\begin{align}
\label{omega2}
 \Omega_2 = \left\{\|\eps\|_2^2/n\leq 1.1\sig^2,~\Lambda_{\max}((\Sig^{n}_{S,S})^{-1})\leq 4/C_{\min}\right\}.
\end{align}
In $\Omega_2$,
\begin{align*}
sgn(\check{b}_S)^T(X_S^TX_S)^{-1}sgn(\check{b}_S)&\leq\|sgn(\check{b}_S)\|_2^2\|(X_S^TX_S)^{-1}\|^2_2\leq \frac{4s}{nC_{\min}}.
\end{align*}
And hence
\begin{align*}
 \Var(T^o_k|X_S,\eps) \leq \Sig_{k,k}\left(\frac{4s}{nC_{\min}}+ \frac{1.1\sig^2}{n\lam^2}\right).
 \end{align*}
Thus, 
\begin{align*}
 & \P\left(\max_{k\in S^c}|T^o_k|\geq \max_{k\in S^c}|\Sig_{k,S}\Sig^{-1}_{S,S}sgn(\check{b}_S)|+t\big|\Omega_2\right) \\
 &\quad\leq 2(p-s)\exp\left\{-\frac{t^2}{2K_1(\frac{4s}{nC_{\min}}+ \frac{1.1\sig^2}{n\lam^2})}\right\}.
\end{align*}
By setting $t = \frac{1-\phi}{2}$, we have for
\[
   n\geq \frac{128K_1s\log p}{C_{\min}(1-\phi)^2}~~\text{and}~\lam\geq \frac{8.8\sig}{1-\phi}\sqrt{\frac{K_1\log p}{n}},
\]
it holds that
\[
   \P\left( \max_{k\in S^c}|T^o_k|>\frac{1+\phi}{2}\big|\Omega_2\right)\leq 2\exp(-\log(p-s)).
 \]
By Lemma \ref{alem2},
\[
  P(\Omega_2)\geq 1-c_1/n-\exp(-c_2n).
\]
Therefore,
\[
  \P(\widehat{S}\subseteq S)\geq 1-c_0/p-c_1/n-\exp(-c_2n).
\]
In the event that $\{\widehat{S}\subseteq S\}$, the KKT condition of (\ref{beta-init}) yields
\begin{equation}
\label{eq3-2}
    \hat{u}_S=(\Sig^{n}_{S,S})^{-1}W^n_{S}-\lam (\Sig^{n}_{S,S})^{-1}sgn(\hat{\beta}_S).
\end{equation}

(ii) Noticing that for any $k\in S$, $\P\{\Sig^n_{k,k}\geq C_{\min}/2\}\geq 1-\exp(-c_2n)$ by Lemma \ref{alem2} and 
\begin{align*}
&\hat{u}_k=\frac{1}{\Sig^n_{k,k}}W^n_k-\frac{\lam}{\Sig^n_{k,k}}sgn(\hat{\beta}_k)-\frac{1}{\Sig^n_{k,k}}\Sig^n_{k,A_k}\hat{u}_{A_k}\\
&=\underbrace{\frac{1}{\Sig^n_{k,k}}W^n_k-\frac{\lam}{\Sig^n_{k,k}}sgn(\hat{\beta}_k)}_{R_{1,k}}-\underbrace{\frac{1}{\Sig^n_{k,k}}\Sig^n_{k,A_k}(\hat{\beta}_{A_k}-\hat{\beta}^{(k)}_{A_k})}_{R_{2,k}}-\underbrace{\frac{1}{\Sig^n_{k,k}}\Sig^n_{k,A_k}(\hat{\beta}^{(k)}_{A_k}-\beta_{A_k})}_{R_{3,k}},
\end{align*}
where $\hat{\beta}^{(k)}_{A_k}$ is the leave-one-out estimator (\ref{beta-loo}).
Using the sub-Gaussian property of $\eps$, it is easy to show that for some large enough $c_1$,
\[
  \P\left(\max_{k\in S} |R_{1,k}|>c_1(\lam+\sqrt{\log s/n})/C_{\min}\right)\leq 2\exp(-c_3\log p).
\]
For $R_{2,k}$, it holds that for some large enough $c_2$,
\[
 \max_{k\in S} |R_{2,k}|\leq \max_{k\in S}\left\|\frac{1}{\Sig^n_{k,k}}\Sig^n_{k,A_k}\right\|_2\|\hat{\beta}_{A_k}-\hat{\beta}^{(k)}_{A_k}\|_2\leq 4\sqrt{\frac{C_{\max}}{C_{\min}}}\rho_S\lam
\]
with probability at least $1-c_1/p-\exp(-c_2n)-c_3/n$ due to Lemma \ref{alem2} and \ref{lem-loo}.
For $R_{3,k}$, it holds that
\begin{align*}
&\max_{k\in S} \left|\Sig^n_{k,A_k}(\hat{\beta}^{(k)}_{A_k}-\beta_{A_k})\right|\\
&\leq \max_{k\in S}\left| \Sig_{k,A_k}\Sig^{-1}_{A_k,A_k}\Sig^n_{A_k,A_k}(\hat{\beta}^{(k)}_{A_k}-\beta_{A_k})\right| +\max_{k\in S} \frac{1}{n}\left|\langle X_k^{\perp},X_{A_k}\rangle(\hat{\beta}^{(k)}_{A_k}-\beta_{A_k})\right|\\
&\leq \rho_S\max_{k\in S}\|\Sig^n_{A_k,A_k}(\hat{\beta}^{(k)}_{A_k}-\beta_{A_k})\|_{\infty} +O_P(\max_{k\in S}\|X_{A_k}(\hat{\beta}^{(k)}_{A_k}-\beta_{A_k})\|_2/n)\\
&\leq \rho_S\lam+\frac{c_5\sqrt{s}\lam}{\sqrt{n}}
\end{align*}
with probability $1-c_2s/n-\exp(-c_3n)$ if $s\log p\leq c_4n$ for a small enough constant $c_4$.
The last step is due to the KKT condition of $\hat{\beta}^{(k)}_{A_k}$ and (\ref{eq-tlem2-5}).

Together we have, for some large enough constant $C$, we have
\[
  \max_{k\in S}|u_k|\leq C\rho_S\lam
\]
with probability at least $1-c_1/p-\exp(-c_2n)-c_3/n$ for some positive constants $c_1-c_3$.

\end{proof}

\subsection{Proof of Lemma \ref{lem1b}}
\label{pf-lem1b}
For $\hat{\beta}^{(k)}$ defined in (\ref{beta-loo}), define a constrained noiseless Lasso estimator 
\begin{equation}
\label{lem2b-eq0}
  \hat{\beta}^{*,(k)}=\argmin_{b\in\R^p}\left\{\frac{1}{2n}\|X_{S}\hat{\beta}^{(k)}_S-X_{S}b\|^2_2+\lam\|b\|_1:b_k=\beta_k, ~b_{S^c}=0\right\}.
\end{equation}
We first prove the following technical lemma.
\begin{lemma}
\label{lem-loo2}
For $\hat{\beta}^{*,(k)}$ defined in (\ref{lem2b-eq0}), it holds that if $n\geq c_0s\log p$ for large enough $c_0$, then with probability at least $1-c_3/p-\exp(-c_4n)-c_5/n$ there exists a large enough constant $c_1$ such that 
\[
 \max_{k\in S}\|\hat{\beta}^{*,(k)}_{S}-\hat{\beta}^{*}_{S}\|_2\leq c_1\rho_S\lam.
\]
\end{lemma}
\begin{proof}[Proof of Lemma \ref{lem-loo2}]
In event $\Omega^*_0$, KKT condition of (\ref{lem2b-eq0}) and that of  $\hat{\beta}^*$ give that
\begin{align*}
-\lam sgn(\hat{\beta}^{*,(k)}_{S})+\lam sgn(\hat{\beta}^{*}_{S})&=\Sig^n_{S,S}(\hat{\beta}^{*,(k)}_{S}-\hat{\beta}^{(k)}_{S})-\Sig^n_{S,S}(\hat{\beta}^*_{S}-\hat{\beta}_{S})\\
&=\Sig^n_{S,S}(\hat{\beta}^{*,(k)}_{S}-\hat{\beta}^{*}_{S})+\Sig^n_{S,S}(\hat{\beta}_{S}-\hat{\beta}^{(k)}_{S}).
\end{align*}
We arrive at
\begin{align*}
 0&\geq -\lam \langle \hat{\beta}^{*,(k)}_{S}-\hat{\beta}^{*}_{S},sgn(\hat{\beta}^{*,(k)}_{S})- sgn(\hat{\beta}^{*}_{S})\rangle\\
 &(\hat{\beta}^{*,(k)}_{S}-\hat{\beta}^{*}_{S})^T \Sig^n_{S,S}(\hat{\beta}^{*,(k)}_{S}-\hat{\beta}^{*}_{S})-(\hat{\beta}^{*,(k)}_{S}-\hat{\beta}^{*}_{S})^T\Sig^n_{S,S}(\hat{\beta}_{S}-\hat{\beta}^{(k)}_{S}),
\end{align*}
where the first step is by the definition of sub-differential.
Therefore,
\begin{align*}
  &\Lambda_{\min}(\Sig^n_{S,S})\|\hat{\beta}^{*,(k)}_{S}-\hat{\beta}^{*}_{S}\|_2^2\leq \|\hat{\beta}^{*,(k)}_{S}-\hat{\beta}^{*}_{S}\|_2\|\Sig^n_{S,S}\|_2\|\hat{\beta}_{S}-\hat{\beta}^{(k)}_{S}\|_2\\
  &\quad\leq  \frac{\Lambda_{\min}(\Sig^n_{S,S})}{2}\|\hat{\beta}^{*,(k)}_{S}-\hat{\beta}^{*}_{S}\|_2^2+ \frac{2\|\Sig^n_{S,S}\|^2_2}{\Lambda_{\min}(\Sig^n_{S,S})}\|\hat{\beta}_{S}-\hat{\beta}^{(k)}_{S}\|_2^2,
\end{align*}
where the last step follows from the Young's inequality.
In the event that $\{C_1\leq \Lambda_{\min}(\Sig^n_{S,S})\leq \Lambda_{\max}(\Sig^n_{S,S})\leq C_2\}$, we arrive at
\[
   \max_{k\in S}\|\hat{\beta}^{*,(k)}_{S}-\hat{\beta}^{*}_{S}\|_2\leq \frac{C_2}{C_1}\max_{k\in S}\|\hat{\beta}_{S}-\hat{\beta}^{(k)}_{S}\|_2\leq c_1\rho_S\lam
\]
for large enough constant $c$ with probability at least $1-c_3/p-\exp(-c_4n)-c_5/n$. The last step is due to Lemma \ref{lem-loo}.

\end{proof}

\begin{proof}[Proof of Lemma \ref{lem1b}]
(i)
Define a restricted Lasso problem with observations $(X,y^*)$.
\begin{equation}
\label{cbs-star}
   \check{b}^*_S = \argmin_{b\in \R^s} \left\{\frac{1}{2n}\|y^*-X_Sb_S\|_2^2+\lam\|b_S\|_1\right\}.
\end{equation}
One can use same arguments in the proof of Lemma \ref{lem1a} to get that 
\[
   \P\{\widehat{S}^*\subseteq S\}\geq 1-c_1/p-c_2/n-\exp(-c_3n).
\]
(ii) The KKT condition for the bootstrapped Lasso gives that
\begin{align*}
    \hat{u}^*_S &= -\lam (\Sig^{n}_{S,S})^{-1}sgn(\hat{\beta}_S^*).
\end{align*}

For any $k\in S$, for $\hat{\beta}^{*,(k)}$ defined in (\ref{lem2b-eq0}), we have
\begin{align*}
\hat{u}^*_k&=-\lam  \frac{1}{\Sig^n_{k,k}}sgn(\hat{\beta}_k^*)- \frac{1}{\Sig^n_{k,k}}\Sig^n_{k,A_k}\hat{u}^*_{A_k}\\
&=\underbrace{ -\lam  \frac{1}{\Sig^n_{k,k}}sgn(\hat{\beta}_k^*)}_{R_{1,k}}-\underbrace{ \frac{1}{\Sig^n_{k,k}}\Sig^n_{k,A_k}[\hat{\beta}_{A_k}^*-\hat{\beta}^{*,(k)}_{A_k}] }_{R_{2,k}} -\underbrace{ \frac{1}{\Sig^n_{k,k}}\Sig^n_{k,A_k}[\hat{\beta}^{*,(k)}_{A_k}-\hat{\beta}^{(k)}_{A_k}]}_{R_{3,k}}\\
&\quad-\underbrace{ \frac{1}{\Sig^n_{k,k}}\Sig^n_{k,A_k}[\hat{\beta}^{(k)}_{A_k}-\hat{\beta}_{A_k}]}_{R_{4,k}}.
\end{align*}
Note that $R_{1,k}=O_P(\lam)$.
For $R_{2,k}$, by Lemma \ref{lem-loo2}  we have with probability at least $1-c_4/p-c_5/n-\exp(-c_6n)$, for some large enough constant $c_7$,
\[
  \max_{k\in S}|R_{2,k}|\leq \|\Sig^{n}_{S,S}\|_2\max_{k\in S}\|\hat{\beta}_{A_k}^*-\hat{\beta}^{*,(k)}_{A_k}\|_2=O_P(\rho_S\lam).
\]
For $R_{3,k}$, by the KKT condition of $\hat{\beta}^{*,(k)}$, we have
\begin{align*}
&\max_{k\in S} |R_{3,k}|\\
&\leq \max_{k\in S} \left|\Sig_{k,A_k}\Sig^{-1}_{A_k,A_k}\Sig^n_{A_k,A_k}[\hat{\beta}^{*,(k)}_{A_k}-\hat{\beta}^{(k)}_{A_k}]\right|+\max_{k\in S} \frac{1}{n}|\langle X_k^{\perp},X_{A_k}\rangle[\hat{\beta}^{*,(k)}_{A_k}-\hat{\beta}^{(k)}_{A_k}]|\\
&\leq \rho_S\max_{k\in S}\|\Sig^n_{A_k,A_k}[\hat{\beta}^{*,(k)}_{A_k}-\hat{\beta}^{(k)}_{A_k}]\|_{\infty}+O_P(\max_{k\in S}\|X_{A_k}(\hat{\beta}^{*,(k)}_{A_k}-\hat{\beta}^{(k)}_{A_k})\|_2/n)\\
&=\rho_S\lam+o_P(n^{-1/2}),
\end{align*}
where the second step is due to $\hat{\beta}^{*,(k)}_{A_k}-\hat{\beta}^{(k)}_{A_k}$ is independent of $X_k$ and the last step is due to a similar argument for (\ref{eq-tlem2-5}).
$\max_{k\in S}|R_{4,k}|$ can be similarly bounded as for $R_{2,k}$ using Lemma \ref{lem-loo}. Hence, we have with probability at least $1-c_7/p-c_8/n-\exp(-c_9n)$, there exists a large enough constant $C^*$ such that
\begin{align*}
\max_{k\in S} |\hat{u}^*_k|\leq C^*\rho_S\lam.
\end{align*}
Therefore $\Omega_0^*$ holds with high probability.
In view of $\Omega_0$ (\ref{omega0}), for $j\in S_+$, 
\begin{align*}
&\hat{\beta}_j = \beta_j + \hat{u}_j \geq  \beta_j - \max_{1\leq j\leq p} |u_j| \geq \beta_j -C\rho_S\lam > C'\rho_S\lam~\text{for}~ \beta_j\geq 0.\\
&\hat{\beta}_j = \beta_j + \hat{u}_j \leq  \beta_j +\max_{1\leq j\leq p} |u_j| \leq \beta_j + C\rho_S\lam < -C'\rho_S\lam ~\text{for} ~\beta_j <0.
\end{align*}
In view of $\Omega_0^*$ (\ref{omega0*}) and repeating previous arguments, we can get $sgn(\hat{\beta}_j)= sgn(\hat{\beta}^*_j)$ for $j\in S_+ $.

\end{proof}


\begin{proof}[Proof of Theorem \ref{thm10}]
We first notice that for $\hat{w}_j$ defined in (\ref{eq-z2}) with $\Theta$ unknown,
\[
   \hat{\beta}^{(mDB)}_j-\beta_j=(e^T_j-\hat{w}_j^TX/n)\hat{u}+\frac{\hat{w}_j^T\eps}{n}.
\]
We have
\begin{align*}
\left|(e^T_j-\hat{w}_j^TX/n)\hat{u}\right|&\leq \left\|\hat{w}_j^TX_{-j}/n\right\|_{\infty}\|\hat{u}\|_1\leq \lam_j\|\hat{u}\|_1=O_P(s\lam\lam_j),
\end{align*}
where the second step is by the KKT condition of (\ref{hkappa}) and the last step is by the standard analysis of the Lasso.
Similarly,
\begin{align*}
 |\hat{b}^{*}_j|&=\left|(e^T_j-\hat{w}_j^TX/n)\hat{u}^*\right|\leq \lam_j\|\hat{u}^*\|_1.
\end{align*}
To obtain an upper bound on $\|\hat{u}^*\|_1$, consider the oracle inequality of $\hat{\beta}^*$:
\begin{align*}
 \frac{1}{2n}\|X\hat{u}^*\|_2^2\leq \lam\|\hat{\beta}\|_1-\lam\|\hat{\beta}^*\|_1&\leq \lam\|\hat{u}^*_S\|_1+\lam\|\hat{\beta}_{S^c}\|_1-\lam\|\hat{\beta}^*_{S^c}\|_1\\
 &\leq \lam\|\hat{u}^*_S\|_1+2\lam\|\hat{\beta}_{S^c}\|_1-\lam\|\hat{u}^*_{S^c}\|_1,
\end{align*}
where we use $\|\hat{u}^*_{S^c}\|_1\leq \|\hat{\beta}^*_{S^c}\|_1+\|\hat{\beta}_{S^c}\|_1$ in the last step. If $2\|\hat{\beta}_{S^c}\|_1\geq \|\hat{u}^*_S\|_1$, then the right hand side gives  $\|\hat{u}^*_{S^c}\|_1\leq 4\|\hat{\beta}_{S^c}\|_1$.
Hence,
\[
   \|\hat{u}^*\|_1\leq \frac{9}{2}\|\hat{\beta}_{S^c}\|_1=\frac{9}{2}\|\hat{u}_{S^c}\|_1\leq \frac{9}{2}\|\hat{u}\|_1=O_P(s\sqrt{\log p/n}).
\]
If $2\|\hat{\beta}_{S^c}\|_1< \|\hat{u}^*_S\|_1$, then we arrive at
\begin{align*}
 \frac{1}{2n}\|X\hat{u}^*\|_2^2\leq \lam\|\hat{\beta}\|_1-\lam\|\hat{\beta}^*\|_1
 &\leq 2\lam\|\hat{u}^*_S\|_1-\lam\|\hat{u}^*_{S^c}\|_1,
\end{align*}
which is the usual oracle inequality. Standard analysis gives $\|\hat{u}^*_S\|_1=O_P(s\sqrt{\log p/n})$.
As a result, 
\begin{align*}
   \hat{\beta}^{(BS-DB)}_j-\beta_j&=\hat{\beta}^{(mDB)}_j-\beta_j-\hat{b}^*_j=\frac{\hat{w}_j^T\eps}{n}+O_P\left(\frac{s\log p}{n}\right).
\end{align*}

\end{proof}

\subsection{Proof of Lemma \ref{lem0b}}
\label{pf-lem0b}
\begin{proof}[Proof of Lemma \ref{lem0b}]

We note that
\[
   P_{\widehat{S}}^{\perp}y=P_{\widehat{S}}^{\perp}(y-X_{\widehat{S}}\hat{\beta}_{\widehat{S}})=P_{\widehat{S}}^{\perp}(y-X\hat{\beta}).
\]
When Lemma \ref{lem1b} holds true,
{\small
\[
   \frac{1}{n-|\widehat{S}|}\|P_{S_+}^{\perp}(y-X\hat{\beta})\|_2^2\leq \hat{\sig}^2=\frac{1}{n-|\widehat{S}|}\|P_{\widehat{S}}^{\perp}(y-X\hat{\beta})\|_2^2\leq\frac{1}{n-|\widehat{S}|}\|P_S^{\perp}(y-X\hat{\beta})\|_2^2.
\]
}
For the left hand side,
\begin{align*}
& \frac{1}{n-|\widehat{S}|}\|P_{S_+}^{\perp}(y-X\hat{\beta})\|_2^2=\frac{1}{n-|\widehat{S}|}\|P_{S_+}^{\perp}\eps\|_2^2 + \frac{1}{n-|\widehat{S}|}\|P_{S_+}^{\perp}X\hat{u}\|_2^2\\
& \quad-\frac{2}{n-|\widehat{S}|}\langle P_{S_+}^{\perp}\eps,P_{S_+}^{\perp}X\hat{u}\rangle,
\end{align*}
where $\hat{u}=\hat{\beta}-\beta$. We have
\begin{align*}
    \frac{1}{n}\|P_{S_+}^{\perp}X\hat{u}\|_2^2&\leq  \frac{1}{n}\|P_{S_+}^{\perp}X_{S_-}\hat{u}_{S_-}\|_2^2\leq \|\Sig^{n}_{S_-,S_-}\|_2 \|\hat{u}_{S_-}\|_2^2 =O_P(s_-\rho^2_S\log p/n).
\end{align*}
We can obtain another bound by noticing $\|\hat{u}_{S_-}\|_2^2\leq \|\hat{u}_S\|_2^2=O_P(s\log p/n).$
Therefore,
\[
   \frac{1}{n}\|P_{S_+}^{\perp}X\hat{u}\|_2^2=O_P\left( \min\{\rho^2_Ss_-,s\}\log p/n\right).
\]
Moreover,
\begin{align*}
   \left|\frac{1}{n}\langle P_{S_+}^{\perp}\eps,P_{S_+}^{\perp}X_{S_-}\hat{u}_{S_-}\rangle\right|&\leq \max_{j\in S_-}\left|\frac{1}{n}\langle P_{S_+}^{\perp}\eps,P_{S_+}^{\perp}X_{j}\rangle\right|\|\hat{u}_{S_-}\|_1=O_P\left(\sqrt{\frac{\log p}{n}}\right)\|\hat{u}_{S_-}\|_1,
\end{align*}
where
\begin{align*}
   &\|\hat{u}_{S_-}\|_1\leq \min\{s_-\|\hat{u}_{S_-}\|_{\infty},\sqrt{s_-}\|\hat{u}_S\|_2\}=O_P\left(\min\{s_-\rho_S\lam,\sqrt{s_-s}\sqrt{\log p/n}\}\right).
\end{align*}
Therefore,
\[
     \left|\frac{1}{n}\langle P_{S_+}^{\perp}\eps,P_{S_+}^{\perp}X_{S_-}\hat{u}_{S_-}\rangle\right|=O_P\left( \min\{\rho_Ss_-,s\}\frac{\log p}{n}\right).
\]
Finally,
\begin{align*}
\frac{1}{n-|\widehat{S}|}\|P_{S_+}^{\perp}\eps\|_2^2 =\sig^2\frac{n-s_+}{n-\widehat{S}}+O_P((n-s)^{-1/2})=\sig^2+O_P(s_-/n+n^{-1/2})
\end{align*}
for $n\gg s\geq |\widehat{S}|$.  Hence,
\[
   \hat{\sig}^2\geq \sig^2-O_P\left(n^{-1/2}+\frac{\min\{\rho_S^2s_-,s\}\log p}{n}\right).
\]
By similar arguments, one can show that
\[
    \frac{1}{n-|\widehat{S}|}\|P_S^{\perp}y\|_2^2=\sig^2+O_P\left(n^{-1/2}+\frac{s_-}{n}\right)
\]
for $n\gg s$. In view of the first inequality in the proof, the proof is complete.
\end{proof}

\section{More simulation results}
We report some numerical results with a ``nonparametric configuration'' of true coefficients. Specifically, $\beta_k=k^{-0.5}$ for $k=1,\dots,s$. We consider the true covariance matrix of $X$ being identity or $\Sig^o$ for $\Sig^o$ defined in Section \ref{sec1-3}. We report the average coverage probabilities on and off the true support and the average confidence interval lengths on and off the true support. In Table \ref{table5}, we see that the proposed BS-DB has average coverage probabilities close to the nominal level when $\Sig=I_p$ and the debiased Lasso method has coverage lower than the nominal level. For $\Sig=\Sig^o$, the coverage probabilities given by BS-DB is closer to the nominal level than those given by the debiased Lasso but both methods have the average coverage probabilities for $\beta_j$, $j\in S$, lower than the nominal level.
\begin{table}[ht]
\centering
\begin{tabular}{cc|cccc|cccc}
\hline
  & & \multicolumn{4}{c|}{$\Sig=I_p$} &  \multicolumn{4}{c}{$\Sig=\Sig^o$} \\
  \hline
  $s$ & j & cov.bsdb & cov.db & se.bsdb & se.db & cov.bsdb & cov.db & se.bsdb & se.db \\ 
  \hline
 4 & $S$ & 0.919 & 0.891 & 0.15 & 0.15 & 0.893 & 0.810 & 0.13 & 0.13 \\ 
   4 & $S^c $   & 0.943 & 0.887 & 0.15 & 0.15   & 0.943 & 0.908 & 0.13 & 0.13 \\ 
   \hline
   8 & $S$ & 0.944 & 0.910 & 0.15 & 0.15 & 0.883 & 0.844 & 0.15 & 0.15\\ 
   8 & $S^c$             & 0.950 & 0.907 & 0.15 & 0.15  & 0.948 & 0.904 & 0.15 & 0.15\\ 
  \hline
  12 & $S$ & 0.939& 0.894 & 0.15 & 0.15 & 0.898 & 0.869 & 0.16 & 0.15 \\ 
   12 & $S^c$   & 0.949 & 0.894 & 0.15 & 0.15   & 0.944 & 0.914 & 0.16 & 0.15\\ 
   \hline
\end{tabular}
\caption{Coverage probabilities of BS-DB (cov.bsdb) and debiased Lasso (cov.db) and standard errors of BS-DB (se.bsdb) and debiased Lasso (se.db). The rows with $j\in S$ are the average results for the coefficients on the true support and the rows with $j\in S^c$ are the average results for the coefficients off the true support.}
\label{table5}
\end{table}

\end{document}